\documentclass{article}
\usepackage{amsmath, amssymb}
\usepackage{amssymb,pb-diagram}
\usepackage{amscd}

\usepackage{fancybox}
\newcommand{\RR}{\mathbb R}

\newcommand{\ZZ}{\mathbb Z}
\newcommand{\PP}{\mathbb P}
\newcommand{\QQ}{\mathbb Q}
\newcommand{\CC}{\mathbb C}

\newcommand{\calE}{{\mathcal E}}

\newcommand{\MW}{\mathop {\rm MW}\nolimits}
\newcommand{\Gal}{\mathop {\rm Gal}\nolimits}

\newcommand{\Hom}{\mathop {\rm Hom}\nolimits}

\newcommand{\NS}{\mathop {\rm NS}\nolimits}

\newcommand{\rank}{\mathop {\rm rank}\nolimits}
\newcommand{\Pic}{\mathop {\rm Pic}\nolimits}

\newcommand{\Supp}{\mathop {\rm Supp}\nolimits}
\newcommand{\Sing}{\mathop {\rm Sing}\nolimits}

\newcommand{\calD}{\mathcal D}

\newtheorem{thm}{Theorem}[section]
\newtheorem{cor}{Corollary}[section]
\newtheorem{prop}{Proposition}[section]
\newtheorem{lem}{Lemma}[section]
\newtheorem{defin}{Definition}[section]
\newtheorem{exmple}{Example}[section]
\newtheorem{rem}{Remark}[section]
\newtheorem{qz}{Question}[section]
\newtheorem{prbm}{Problem}[section]
\newcommand{\qed}{\hfill $\Box$}

\newcommand{\proof}{\noindent{\textsl {Proof}.}\hskip 3pt}
\newcommand{\proofend}{\qed \par\smallskip\noindent}

\renewcommand{\thesubparagraph}{\theparagraph.\@arabic\c@subparagraph}
\begin{document}
  
  \begin{center}
  
  \Large Geometry of irreducible plane quartics \\ and \\ their quadratic residue conics

\bigskip

%\normalsize Dedicated to Professor Du Plessis on his sixties birthday.

\bigskip
\large 
Hiro-o TOKUNAGA

\end{center}
\normalsize
%\titlerunning{ }  % abbreviated title (for running head)
%\authorrunning{ }  % abbreviated title (for running head)
%\pagestyle{headings}  % switches on printing of running heads
%%\thispagestyle{empty} %%nopage

%\institute{
%Department of Mathematics and Information Sciences\\
%Graduate School of Science and Engineering,\\
%Tokyo Metropolitan University\\
%1-1 Minami-Ohsawa, Hachiohji-shi 192-0397 JAPAN \\
%{\tt tokunaga@tmu.ac.jp, kameda-yuuiti@ed.tmu.ac.jp} \newline
%\and
%Graduate School of Informatics, Kyoto University\\
%Yoshida Honmachi, Sakyo-ku, Kyoto 606-850 JAPAN\\
%{\tt akihiro@i.kyoto-u.ac.jp}
%}

%\maketitle

\begin{abstract}  
Let $D$ be an irreducible plane curve in $\PP^2$. In this article, we first introduce a notion of
a quadratic residue curve $\bmod D$, and study quadratic residue concis $C$ $\bmod$ an irreducible
quartic curve $Q$. As an application, we study a dihedral cover of $\PP^2$ with branch locus
$C+Q$ and give two examples of Zariski pairs as by-products. 
\end{abstract}

  \large
  
  {\bf Introduction}
 
  \normalsize
 In this article, we study geometry of irreducible plane quartic $Q$ and a smooth conic $C$ which
 is tangent to $Q$ with even order at each point in $C\cap Q$.  Geometry of a smooth plane quartic
 and its bitangent lines is a classical object and well-studied by many mathematicians from
  various points of view. We hope that this article adds another interesting topic to geometry of
  plane quartics.
  All varieties throughout this paper are defined over the field of complex numbers,  ${\mathbb C}$.
In order to explain our motivation and results on the above subject, let us start with introducing some
notion and definition. 

   Let $\Sigma$ be a smooth projective surface.
  Let   $f' : Z' \to \Sigma$ be a double cover of $\Sigma$, i.e.,  $Z'$ is a normal surface and 
   $f'$ is a finite surjective morphism of degree $2$.  We denote its canonical resolution by
   $\mu : Z \to Z'$ (see \cite{horikawa} for the canonical resolution). Note that
   that $\mu$ is the identity if $Z'$ is smooth.
   We put $f := f'\circ\mu$. We denote the involution on $Z$ induced by
   the covering transform of $f'$ by $\sigma_f$. The branch locus $\Delta_{f'}$ of $f'$ is 
   the subset of $\Sigma$ consisting of points $x$ such that
   $f'$ is not locally isomorphic over $x$. Similarly we define the 
   branch locus $\Delta_f$ of  $f$. Note that $\Delta_{f'} = \Delta_{f}$.

   \begin{defin} \label{def:splitting}{\rm Let $D$ be an irreducible curve on $\Sigma$.
   We call $D$ a splitting curve with respect to  $f$ if $f^*D$ is of the 
   form
   \[
   f^*D = D^+ + D^- + E,
   \]
   where $D^+ \neq D^-$, $\sigma_f^* D^+  = D^-$, $f(D^+) = f(D^-) = D$ and
   $\Supp(E)$ is contained in the exceptional set of $\mu$.  If the double cover
   $f : Z \to \Sigma$ is determined by its branch locus $\Delta_f$, i.e., any double
   cover with branch locus $\Delta_f$ is isomorphic to $Z'$ over
   $\Sigma$,  and $D$
   is a splitting curve with respect to $f$, we say that
   ^^ ^^ {\it $\Delta_f$ is a quadratic residue curve mod $D$}".
   }
   \end{defin}

   \begin{rem}\label{rem:no-1}{\rm
   \begin{itemize}
  \item Note that if $\Sigma$ is simply connected, then any double cover of
    $\Sigma$ is determined by its branch locus.
    
    \item In our previous results on  dihedral covers and their application to the study of 
  the topology of the complements of plane curves, we see that splitting curves play important roles and that it
  is indispensable to  know the property of them. (see \cite{artal-tokunaga}, \cite{tokunaga94}, \cite{tokunaga05}, for example). This is our first motivation to study splitting curves.

  \item Our  terminology comes from elementary number theory.  Let $m$ be
   a square free positive integer,  let $p$ be an odd prime with $p \not| m$ and let 
   ${\mathcal O}_{{\mathbb Q}(\sqrt{m})}$ be the integer ring
   of ${\mathbb Q}(\sqrt{m})$.  It is well-know that the ideal $(p)$ generated by
   $p$ in ${\mathcal O}_{{\mathbb Q}(\sqrt{m})}$ satisfies the following properties
   (See \cite[Proposition 13.1.3]{ireland-rosen}, p.190, for example):
   
   \begin{itemize}
   \item If $m$ is a quadratic residue mod $p$, then
   $(p) = {\mathfrak p}_1{\mathfrak p}_2$, where ${\mathfrak p}_i$ $(i = 1, 2)$ are
    distinct prime ideals.
    
    \item If $m$ is not a quadratic residue mod $p$, then
    $(p)$ is a prime ideal.   
  \end{itemize}
  \end{itemize}
        }
    \end{rem}

   Suppose that $f : Z \to \Sigma$ is uniquely
  determined by $\Delta_f$. 
  Likewise the Legendre symbol in elementary number theory, we here introduce a notation to describe if $\Delta_f$ is
  a quadratic residue mod $D$ or not. For an irreducible curve $D$ on $\Sigma$, 
  we put
  \[
  (\Delta_f/D) = \left \{ \begin{array}{cc}
                         1 & \mbox{if $\Delta_f$ is a quadratic residue curve mod $D$} \\
                         -1 &\mbox{if $\Delta_f$ is not a quadratic residue curve mod $D$}
                                \end{array}
                                \right .
   \]

  As $\PP^2$ is simply connected, any double cover of $\PP^2$ is just determined by its branch locus.
 On the other hand, any reduced plane curve $B$  of even degree can be the branch locus of
 a double cover. Hence for any irreducible plane curve $D$, one can consider $(B/D)$. 
 
In this article, we consider the case when \textit{any point $x \in B\cap D$ is a smooth point of both $B$ and $D$}.
 For such a case, if the intersection multiplicity at some point in $B\cap D$ is odd, then we infer that
 $(B/D) = -1$. This leads us to introduce  a notion of \textit{even tangential curve}. 
 
  \begin{defin}\label{def:etc}{\rm Let $D_1$ and $D_2$ are  reduced divisors on 
 a smooth projective surface without any common 
 irreducible component.
 We say that $D_1$ and $D_2$ are  even tangential or $D_1$ (resp. $D_2$) is
 even tangential to $D_2$ (resp. $D_1$) if
 \begin{enumerate}
 
 \item[(i)] For $\forall P \in D_1\cap D_2$, $P \not\in \Sing(D_1)\cup\Sing(D_2)$, and
 
 \item[(ii)]  the intersection multiplicity of $D_1$ and $D_2$ at $P$, $I_P(D_1, D_2)$,
  is even for $\forall P \in D_1\cap D_2$. 
 \end{enumerate}
 Note that we do not pay attention to $\sharp (D_1\cap D_2)$ to define even tangential
 curves.
 }
 \end{defin}
 
 Now our basic problem can be formulated as follows:

 \begin{prbm}\label{pbm:general}{\rm Let $B$ be a reduced plane curve of even degree.
 
   \begin{enumerate}
  
  \item[(i)] Find an even tangential curve $D$ to $B$ and  determine the value of $(B/D)$.
  
  \item[(ii)] What can we say about the topology of $\PP^2\setminus(B + D)$ from the 
  value of $(B/D)$?
  
  \end{enumerate}
  }
  \end{prbm}
   
  As a first step, we consider the case when $B$ is a smooth conic $C$.  Suppose that  $D$ is an irreducible plane curve which is even tangential to $C$.  We easily see the
  following:
  
  \begin{itemize}
  
  \item If  $\deg D = 1, 2$, we have $(C/D) = 1$.
  
   \item If  $\deg D = 3$, we have 
       \begin{itemize}
         \item $(C/D) = -1$ if $D$ is smooth, and
         \item $(C/D) = 1$ if $D$ is a nodal cubic.
         \end{itemize}
         Note that there is no even tangential cuspidal cubic to $C$.
  \end{itemize}
  
  Hence the case of $\deg D = 4$ seems to be the first interesting case. Now let us restate
  our exact problems which we consider in this article:
  
  \begin{prbm}\label{pbm:specific}{\rm Fix an irreducible quartic  $Q$. 
  
  \begin{enumerate}
  
   \item[(i)]  Find even tangential conics $C$ to $Q$ and determine the value $(C/Q)$.
   
   \item[(ii)] Does the value $(C/Q)$ affect the topology of $\PP^2\setminus (C+Q)$?
   \end{enumerate}
  }
  \end{prbm}
 
 In this article, we first consider Problem~\ref{pbm:specific} (i) and give a formula to determine
 $(C/Q)$ (see Theorem~\ref{thm:qrc}).
We next count the number of even tangential conics passing through a smooth point $x$ on $Q$. 
 Now our result is as follows:
 
 \begin{thm}\label{thm:main}{ Choose a smooth point $x$ of $Q$. There exist
 finitely many (possibly no) even tangential conics $C$ through $x$ and we have the
 following table:

\begin{itemize}

\item $\Xi_Q$: the set of types of singularities of $Q$. Note that $Q$ has at worst simple singularities
and we use the notation in \cite{bpv}  in order to describe the type of a singularity.

\item $l_x\cap Q$: This shows how $l_x$ meets $Q$.  We use the following notation to describe it.

      \begin{itemize}
      \item $s$: $I_x(l_x, Q) = 2$ or $3$, and $l_x$ meets $Q$ transversely at other point(s).

      \item $b$: $l_x$ is either bitangent line through $x$ or $I_x(l_x, Q) = 4$.

      \item $sb$: $I_x(l_x, Q) = 2$ and $l_x$ passes through a double point of $Q$.
      \end{itemize}

\item  {\rm ETC}: the set of even tangential conics passing through $x$ and $\sharp \mbox{\rm ETC}$ denotes its cardinality.

\item {\rm QRETC}:  the set of even tangential conics passing through $x$ with $(C/Q) = 1$ and
$\sharp \mbox{\rm QRETC}$ denotes its cardinality.

\item We omit cases of $(\Xi_Q$,  $l_x\cap Q)$ which do not occur.  For example, the case of 
$(\Xi_Q, l_x\cap Q) = (A_6, b)$ is omitted,  as such a case does not occur.

\end{itemize}

\rm
\begin{center}
\begin{tabular}{|c|c|c|c|c|} \hline 
 No. & $\Xi_Q$ & $l_x\cap Q$ & $\sharp \mbox{\rm ETC}$ & $\sharp \mbox{\rm QRETC}$ \\ \hline
 1  & $A_6$    &     $s$   & $0$ & $0$    \\ 
 2 & $A_6$    &      $sb$  & $0$ & $0$ \\ \hline
 3 & $E_6$    &      $s$    & $0$ & $0$ \\
 4 & $E_6$    &      $b$   & $0$ & $0$   \\ \hline    
 5 & $A_5$   &       $s$   & $1$ & $1$  \\
 6 & $A_5$    &      $b$   & $1$ & $1$  \\
 7 & $A_5$    &      $sb$  & $0$ & $0$   \\ \hline
 8 & $D_5$    &      $s$    & $1$ & $1$  \\
   9 & $D_5$    &      $b$    & $0$ & $0$  \\ \hline
    \end{tabular}
 
 \begin{tabular}{|c|c|c|c|c|} \hline 
 No. & $\Xi_Q$ & $l_x\cap Q$ & $\sharp \mbox{\rm ETC}$ & $\sharp \mbox{\rm QRETC}$ \\ \hline

 10 & $D_4$    &      $s$    & $3$ & $3$  \\
 11 & $D_4$    &      $b$    & $0$ & $0$ \\ \hline
 12 & $A_4 + A_2$ &$s$ &$0$ & $0$  \\
  13 & $A_4 + A_2$ &$sb$ &$0$ & $0$ \\ \hline
  14 & $A_4 + A_1$& $s$ & $0$ & $0$ \\
  15 & $A_4 + A_1$ &   $b$ & $0$ & $0$  \\
16 & $A_4 + A_1$ & $sb$ & $0$ & $0$  \\
17 & $A_4 + A_1$&  $sb$ & $0$ & $0$  \\ \hline    
18 & $A_3 + A_2$ & $s$ &$1$ & $1$ \\   
19 & $A_3 + A_2$ & $sb$ & $0$ & $0$   \\
20 & $A_3 + A_2$ & $sb$ & $1$ & $1$  \\ \hline 
21 & $A_3 + A_1$ & $s$ & $2$ & $2$  \\                                                              
22 & $A_3 + A_1$ & $b$ & $1$ & $1$ \\   
23 & $A_3 + A_1$ & $sb$ & $1$ & $1$  \\   
24  & $A_3 + A_1$ & $sb$ & $0$ & $0$  \\   \hline
25 & $3A_2$ & $s$ & $0$ & $0$  \\
26 & $3A_2$ & $b$ & $0$ & $0$  \\ \hline
27 & $2A_2 + A_1$ & $s$ & $0$ & $0$  \\
28 & $2A_2 + A_1$ & $b$ & $0$ & $0$ \\
29 & $2A_2 + A_1$ & $sb$ & $0$ & $0$  \\ \hline
30 & $A_2 + 2A_1$ & $s$ & $1$ & $1$ \\ 
 31& $A_2 + 2A_1$ & $b$ & $0$ & $0$   \\
 32 & $A_2 + 2A_1$ & $sb$ & $0$ & $0$ \\ 
 33 & $A_2 + 2A_1$ & $sb$ & $1$ & $1$  \\ \hline      
 34 & $3A_1$ & $s$ & $4$ & $4$ \\
 35 & $3A_1$ & $b$ & $1$ &   $1$  \\                                                                                                                                                                                                                
36 & $3A_1$ & $sb$ & $2$ & $2$  \\ \hline                                                               
37 & $A_4$ & $s$ &$3$ & $0$ \\ 
 38 & $A_4$ & $b$ &$1$ &  $0$ \\                                                            
39 &  $A_4$ & $sb$ & $1$  &  $0$  \\ \hline
40 & $A_3$ &  $s$ & $7$ & $1$ \\   
41 & $A_3$ & $b$ & $2$ &  $0$ \\ 
42 & $A_3$ & $sb$ & $4$ & $1$  \\ \hline
43 & $2A_2$ & $s$ & $3$ & $0$ \\
44 & $2A_2$ & $b$ & $3$ & $0$  \\
45 & $2A_2$ & $sb$ & $1$ &    $0$  \\ \hline                                                                       
46 & $A_2 + A_1$ & $s$ & $6$ & $0$  \\ 
 47 & $A_2 + A_1$ & $b$ & $3$ & $0$  \\
 48 & $A_2 + A_1$ & $sb$ & $3$ & $0$  \\ 
  49 &  $A_2 + A_1$ & $sb$ & $2$ & $0$  \\ \hline                                                                          
 50 & $2A_1$ & $s$ & $13$ & $1$ \\                                                                       
51 & $2A_1$ & $b$  & $6$ & $0$  \\ 
52 & $2A_1$ & $sb$ & $7$ & $1$ \\ \hline  
\end{tabular}

\begin{tabular}{|c|c|c|c|c|} \hline 
 No. & $\Xi_Q$ & $l_x\cap Q$ & $\sharp \mbox{\rm ETC}$ & $\sharp \mbox{\rm QRETC}$ \\ \hline
53 & $A_2$ & $s$ & $15$ & $0$ \\
54 & $A_2$ & $b$ & $6$ & $0$ \\
55 & $A_2$ & $sb$ & $10$ & $0$ \\ \hline
56 & $A_1$ & $s$ & $30$ & $0$ \\
57 & $A_1$ & $b$ & $15$ & $0$ \\
58 & $A_1$ & $sb$ & $20$ & $0$ \\ \hline   
59 & $\emptyset$ & $s$ &$63$ & $0$ \\
60 & $\emptyset$ & $b$ & $36$ & $0$ \\ \hline                                                                                                                                       
\end{tabular}                                                                                             
\end{center}
 }
 \end{thm}
 
 Note that there exist both quadratic and non-quadratic residue even tangential conics
 to $Q$ for the cases 40, 42, 50 and 52.  These cases are interesting when
 we consider 
 Problem~\ref{pbm:specific} (ii). In fact, we study dihedral covers of $\PP^2$ whose
 branch loci are $C +  Q$, and have the following result(
 see \S 3  for the notations on dihedral covers):
 
 \begin{thm}\label{thm:main2} {Let $Q$ be an irreducible quartic,  let $C$ be an even
 tangential conic to $Q$ and let $f_C : Z_C \to \PP^2$ be a double cover with $\Delta_{f_C} = C$. If there exists a ${\mathcal D}_{2p}$-cover $\pi: S \to \PP^2$ with
$\Delta_{\pi} = C + Q$ for an odd prime $p \ge 5$, then  we have the following:

\par\medskip

$(i)$ $D(X/\PP^2) = Z_C \cong \PP^1\times\PP^1$, i.e., $\pi$ is branched at $2C + pQ$.

\par\medskip

$(ii)$ $(C/Q) = 1$. 
Moreover, if we put  $f_C^*Q = Q^+ + Q^-$, then $Q^+ \sim Q^- \sim (2,2)$.

Conversely, if the second condition holds, then there exist ${\mathcal D}_{2n}$-covers $\pi_n: 
S_n \to \PP^2$ branched at $2C + nQ$ for any $n \ge 3$.
}
\end{thm}
 From Theorem~\ref{thm:main2}, we have the following corollary:

\begin{cor}\label{cor:pi1}{There exists an epimorphism 
from the fundamental group $\pi_1(\PP^2\setminus (C+Q), \ast)$ to ${\mathcal D}_{2p}$
for an odd prime $p \ge 5$ if and only if $(C/Q) = 1$ and $Q^+ \sim Q^-$
}
\end{cor}

 This paper consists of $5$ sections. In \S 1,  we start with preliminaries on theory of elliptic surface.  We prove Theroem~\ref{thm:main} in \S 2.
 In \S 3, we give a summary on branched Galois covers, mainly dihedral covers. We prove
 Theorem~\ref{thm:main2} in \S 4.
 In \S 5, we consider an application of Theorem~\ref{thm:main2} and give two examples of 
 Zariski pairs.

\bigskip

\textbf{Acknowledgement.} The most of part of this article was done during the author's visit to
Ruhr Universi\"at Bochum under the support of SFB/TR 12. The author thanks  Professor A. Huckleberry for his arrangement and hospitality. The author also thanks
the organizers of the symposium ^^ ^^ Singularities in Aarhus" for giving  the author an opportunity to give a talk on 
the subject in this article. 
 
%  \input sec-coqua1.tex
%sec-coqua1.tex

\section{Preliminaries on elliptic surfaces}
\newcommand{\Corr}{\mathop {\rm Corr}\nolimits}

\subsection{Elliptic surfaces.}

We review some general facts on elliptic surfaces. For details, we refer to
\cite{kodaira}, \cite{miranda} and \cite{shioda90}.
Let $\varphi: {\mathcal E} \to C$ be an elliptic surface over a smooth projective curve $C$ with a section
$O$. Throughout this article, we always assume that
\begin{enumerate}
\item[(i)] $\varphi$ is relatively minimal and 
\item[(ii)] there exists at least one singular fiber.
\end{enumerate}

Let $\NS({\mathcal E})$ be the N\'eron-Severi group of ${\mathcal E}$ and let $T_{\varphi}$ be
the subgroup of $\NS(\calE)$ generated by $O$ and all the irreducible components of fibers
of $\varphi$. $T_{\varphi}$ has a canonical basis as follows:

$O$, a general fiber ${\mathfrak f}$, and $\{\Theta_{v,1}, \ldots, \Theta_{v,m_v-1}\}_{v \in R_{\varphi}}$, where

\begin{itemize}

\item $R_{\varphi}:= \{v \in C \mid \mbox{$\varphi^{-1}(v)$ is reducible.}\}$, and

\item we label the irreducible components of $\varphi^{-1}(v)$ as follows:
$\Theta_{v, 0}, \Theta_{v,1}, \ldots, \Theta_{v, m_v-1}$, $\Theta_{v, 0}O = 1$.

\end{itemize}

Let $\MW(\calE)$ be the Mordell-Weil group, the group of sections, of $\calE$,
$O$ being the zero sections. Under these circumstances, we have 

\begin{thm}\label{thm:shioda}{{\rm \cite[Theorem 1.3]{shioda90}} There is 
a natural isomorphism
\[
\MW(\calE) \cong \NS(\calE)/T_{\varphi}.
\]
}
\end{thm}
Also in \cite{shioda90}, a symmetric bilinear form $\langle \, , \, \rangle$,
called the height pairing, on $\MW(\calE)$ is defined by using the intersection
pairing as follows:

For any $s \in \MW(\calE)$, $\langle s, s \rangle \ge 0$ and $=0$ if and only if
$s$ is a torsion. More explicitly, for $s_1, s_2 \in \MW(\calE)$, $\langle s_1, s_2\rangle$ is given by
\[
\langle s_1, s_2 \rangle = \chi({\mathcal O}_{\calE}) + s_1O +s_2O -s_1s_2 - \sum_{v \in R_{\varphi}}
\mbox{Corr}_v(s_1, s_2),
\]
where $\mbox{Corr}_v(s_1, s_2)$ is given by 
\[
\Corr_v(s_1, s_2) = (s_2\Theta_{v,1}, \ldots, s_2\Theta_{v, m_{v-1}})
(-A_v^{-1})\left ( \begin{array}{c}
                            s_1\Theta_{v,1} \\
                             \cdot \\
                             s_1\Theta_{v, m_v-1}
                             \end{array} \right ) ,
\]
and $A_v$ is the intersection matrix $(\Theta_{v,i}\Theta_{v,j})$ ($1\le i, j \le m_v-1$). As for 
explicit values of $\Corr_v(s_1, s_2)$, see Table 8.16 in \cite{shioda90}. 
%In particular, if $\calE$ is a rational elliptic surface, we have a formula
%\[
%\langle s, s \rangle = 2 + 2sO  - \sum_{v \in R_{\varphi}}
%\mbox{Corr}_v(s, s).
%\]
%Here we call $\calE$ a rational elliptic surface if $\calE$ is a rational surface.

\subsection{A ^^ ^^ reciprocity" between sections and trisections on rational ruled surfaces.}

Let $\Sigma_d$ be the Hirzebruch surface of degree $d$ ($d$: even positive
integer). We denote its section with self-intersection number $-d$ and its
fiber of the ruling by $\Delta_{0,d}$ and $F_d$, respectively. Let $\Gamma_d$
be an irreducible curve on $\Sigma_d$ such that

\begin{enumerate}

\item $\Gamma_d \sim 3(\Delta_{0,d} + dF_d)$ and

\item $\Gamma_d$ has at worst simple singularities.

\end{enumerate}

Let $\Delta$ be a section on $\Sigma_d$ such that $(i)$ $\Delta \sim
\Delta_{0,d} + dF_d$ and $(ii)$ $\Delta$ and $\Gamma_d$ are even tangential.

Let $p'_d : S'_d \to \Sigma_d$ be the double cover with branch locus 
$\Delta_0 + \Gamma_d$ and $\mu_d :  S_d \to S'_d$ be the canonical 
resolution and put $p_d := p'_d\circ\mu_d$. Since $\Delta_0 + \Gamma_d$ meets
distinct $4$ points with a general fiber $F_d$, one can easily see that 
the ruling on $\Sigma_d$ induces an elliptic fibration $\varphi_d : S_d \to \PP^1$.
Moreover, by its construction, we infer that 

\begin{enumerate}

\item[(a)] $\varphi_d$ is relatively minimal,

\item[(b)] the preimage of $\Delta_0$ gives a section which we denote by $O$, and

\item[(c)] $\Delta$ gives rise to two sections $s_{\Delta}^+$ and $s_{\Delta}^-$ of 
$\varphi_d$.

\end{enumerate}

Let $\MW(S_d)$ be the group of sections , the Mordell-Weil group, of ${\varphi}_d$, where
$O$ is the zero element.
Let  $q_d : W_d \to \Sigma_d$ be a double cover with branch locus $\Delta_0 + \Delta$.
Note that $q_d$ is uniquely determined by $\Delta_0 + \Delta$ as $\Sigma_d$ is simply connected and that $W_d \cong \Sigma_{d/2}$. Then we have

\begin{thm}\label{thm:recipro}{
\[
\left ( (\Delta_0 + \Delta)/\Gamma_d)\right ) = (-1)^{\varepsilon(s^+_{\Delta})}
\]
where, for a section $s \in \MW(S_d)$, $\varepsilon(s)$ is defined as follows:
\[
\varepsilon(s) = \left \{ \begin{array}{cc}
                                          0 & \mbox{$\exists s_o \in \MW(S_d)$ such that $s = 2s_o$} \\
                                          1 & \mbox{$\not\exists s_o \in \MW(S_d)$ such that $s = 2s_o$} 
                                          \end{array} \right . 
\]
Note that $\varepsilon(s_{\Delta}^+) = \varepsilon(s_{\Delta}^-)$ as 
$s_{\Delta}^+ = - s_{\Delta}^-$ on $\MW(S_d)$. 
}
\end{thm}

\proof  It is enough to show 
\[
\left( (\Delta_0 + \Delta)/\Gamma_d)\right ) = 1 \Leftrightarrow s_{\Delta}^{\pm} \in
2\MW(S_d).
\]
($\Rightarrow$) As we have seen, $W_d \cong \Sigma_{d/2}$. We choose 
 affine open subsets
   $V \subset W_d (\cong \Sigma_{d/2})$, and $U \subset \Sigma_{d}$ as follows:
   
   \begin{enumerate}
   
   \item[(i)] Both $U$ and $V$ are $\CC^2$.
   
   \item[(ii)] We choose affine coordinates $(t,u)$ and $(\tilde{t}, \zeta)$  of $U$ and $V$, respectively,
   in such a way that $q_d$ is given by
      \[
   q_d : (\tilde{t}, \zeta) \mapsto (t, u) = (\tilde{t}, \zeta^2 + f(t)),
   \]
   where $f(t)$ is a polynomial of degree $\le d$. Note that
    with respect to these coordinate $(t, u)$ and $(\tilde t, \zeta)$, $
  \Delta\cap U : u - f(t) = 0$, $\Delta_{0,d}$ corresponds to the section given by $u = \infty$ and 
  the involution $\sigma_{q_d}$ is given by $(\tilde t, \zeta) \mapsto (\tilde t, - \zeta)$.  
  \end{enumerate}
  
  Since $((\Delta_{0,d} + \Delta)/\Gamma_d) = 1$, $q_d^*\Gamma_d$ is of the form $\Gamma^+ +
  \Gamma^-$. Since $\sigma_{q_d}^*\Gamma^+ = \Gamma^-, \sigma_{q_d}^*\Delta_{0, d/2} = \Delta_{0, d/2}$ and $\sigma_{q_d}^*F_{d/2} = F_{d/2}$,  $\Gamma^+ \sim \Gamma^- \sim 3(\Delta_{0,d/2} + d/2 F_{d/2})$.
  Hence we may assume 
  \begin{eqnarray*}
  \Gamma^+ &:& F(\tilde t, \zeta) = \zeta^3 + a_1(\tilde t)\zeta^2 + a_2(\tilde t)\zeta + a_3(\tilde t) =  0 \\
  \Gamma^- &:& -F(\tilde t, -\zeta) = \zeta^3 - a_1(\tilde t)\zeta^2 + a_2(\tilde t)\zeta - a_3(\tilde t)  =  0,
  \end{eqnarray*}
  where $\deg a_k(\tilde t) \le kd/2$. Since $\zeta^2 = u - f(t), t = \tilde t$, we have
  \[
  F(\tilde t, \zeta) = (a_1(t) u - a_1(t) f(t) + a_3(t)) + (u - f(t) + a_2(t))\zeta.
  \]
  As $q_d^*\Gamma = \Gamma^+ + \Gamma^-$, we may assume that $\Gamma_d$ is given by
  \[
  -F(\tilde t, \zeta)F(\tilde t, -\zeta) = (a_1(t) u - a_1(t) f + a_3(t))^2 - (u - f(t) + a_2(t))^2(u - f(t)) = 0.
  \]
  On the other hand, over $U$ is $S'_d$ is given by
  \[
  S'_d|_{{p'_d}^{-1}} : y^2 = (a_1(t) u - a_1(t) f + a_3(t))^2 - (u - f(t) + a_2(t))^2(u - f(t)),
  \]
  and the above equation considered as a Weierstrass equation of the generic fiber, $S_{d,\eta}$,
  of $\varphi_d$. By our construction, $s_{\Delta}^{\pm}$ is given by
  \[
  s_{\Delta}^{\pm}: (f(t), \pm a_3(t)).
  \]
  Put 
  \[
  s_o^{\pm}: (\mp(f(t) - a_2(t)), \pm(a_1(t)a_2(t) - a_3(t)).
  \]
  Then $s_o^{\pm} \in \MW(S_d)$ and by the definition of the group law, we have
  \[
  2s_o^{\pm} = s_{\Delta}^{\pm}.
  \]
  
  ($\Leftarrow$)  We use the affine open sets of $\Sigma_d, W_d$ as before. Suppose that
  $\Gamma_d$ is given by
  \[
  \Gamma_d : F_{\Gamma_d}(t, u) = u^3 + c_1(t)u^2 + c_2(t)u + c_3(t) = 0
  \]
  where $c_k(t) (i = 1, 2, 3)$ are polynomials of degrees $\le kd$. Then $S'_d$ over $U$ is given by
  $y^2 = F_{\Gamma_d}(t, u)$ and,
  as we have seen, this equation can be regarded as a Weierstrass equation of 
  the generic fiber $S_{d, \eta}$.  Since $s_{\Delta}^+O = 0$ and $p_d(s_{\Delta}^+) = \Delta$,
  $s_{\Delta}^+ \in \MW(S_d)$ is given by
  \[
  s_{\Delta}^+: (u, y) = (f(t), g(t)),
  \]
  where $g(t)$ is a polynomial of degree $\le 3d/2$.  Let $s_o \in \MW(S_d)$ such that $2s_o = s_{\Delta}^+$.
  Since $s_o$ is a $\CC(\PP^1) (= \CC(t))$-rational point of $S_{d, \eta}$, there exist $f_o(t), g_o(t) 
  \in \CC(t)$ such that
  \[
  s_o: (u, y) = (f_o(t), g_o(t)).
  \]
  Since $s_{\Delta}^+ O = 0$, by \cite[Theorem 9.1]{kodaira}, we infer that $s_oO = 0$. Therefore
  $f_o(t), g_o(t) \in \CC[t]$ and $\deg f_o \le d, \deg g_o \le 3d/2$. 
 Now let 
  \[
  y = \alpha(t) u + \beta(t), \alpha(t), \beta(t) \in \CC(t)
  \]
  be the tangent line of the elliptic curve $S_{d, \eta}$ over $\CC(t)$ at $s_o$. By the definition of 
  the group law on $S_{d, \eta}$, we have
  \[
  F(t, u) = (\alpha(t)u + \beta(t))^2 + (u- f_o(t))^2(u - f(t)).
  \]
  As $F(t, u), f, f_o \in \CC[t, u]$, we infer that $\alpha(t), \beta(t) \in \CC[t]$. Thus we may assume
  that $\Gamma_d\cap U$ is given by
  \[
   (\alpha(t)u + \beta(t))^2 + (u- f_o(t))^2(u - f(t)) = 0.
   \]
   As $q_d^*\Gamma_d$ on $V$ is given by
   \begin{eqnarray*}
   && (\alpha(t)u + \beta(t))^2 + (u- f_o(t))^2\zeta^2 \\
   &= & \{ (\alpha(t)u + \beta(t)) +\sqrt{-1} (u- f_o(t))\zeta\}\times\{(\alpha(t)u + \beta(t))- \sqrt{-1} (u- f_o(t))\zeta\},
   \end{eqnarray*}
   $\Gamma_d$ is splitting with respect to $q_d$, i.e., $((\Delta_0 + \Delta)/\Gamma_d) = 1$.
   \proofend

\begin{rem}{\rm Theorem~\ref{thm:recipro} can be generalized to the case when $S_d$ has a hyperelliptic
fibration under some restriction. See \cite{tokunaga09}.
}
\end{rem}

\subsection{Double covers of $\PP^2$ branched along quartics and rational elliptic surfaces}

We call an elliptic surface $\calE$ rational, if $\calE$ is a rational surface. Hence 
its base curve of the elliptic fibration  is $\PP^1$. 
Likewise we have done in \cite{tokunaga94}, we associate a reduced plane curve $Q$ in $\PP^2$
and a distinguished smooth point $x$ on $Q$ with a rational elliptic surface $\calE^Q_x$ as follows:

\medskip

Let $\nu_1: \PP_x^2 \to \PP^2$ be a blowing-up at $x$. We denote the proper transform of 
the tangent line $l_x$ at $x$ by $\overline{l}_{x,1}$, and the exceptional curve of $\nu_1$ by $E_{x,1}$.
We next consider  another blowing up $\nu_2: \widehat{\PP}^2 \to \PP^2_x$ at $\overline{l}_{x,1}
\cap E_{x,1}$, and denote the proper transforms of $\overline{l}_{x,1},  E_{x,1}$ and and the
exceptional curve of $\nu_2$ by $\overline{l}_x$, $\overline{E}_{x,1}$, 
 and $E_{x,2}$,
respectively. Let $f' : \calE' \to \widehat{\PP}^2$ be a double cover with branch locus
$\overline{E}_{x,1}$ and $\overline {Q}$, where $\overline{Q}$ is the proper transform
of $Q$ with respect to $\nu_2\circ\nu_1$. Let $\mu_x^Q : \calE_x^Q \to \calE'$ be the canonical
resolution of $\calE'$ and put $f_x^Q:= f'\circ\mu_x^Q$. Then we see that  $\calE_x^Q$ satisfies the following properties:

\begin{enumerate}
\item[(i)] The pencil $\Lambda_x$ of lines through $x$ induces a relatively minimal
elliptic fibration $\varphi_x^Q : \calE_x^Q \to \PP^1$.

\item[(ii)] The preimage of $\overline{E}_{x,1}$ gives rise to a section $O$ of
$\varphi_x^Q$, and the generic fiber has a group structure, $O$ being the zero element.
 Moreover the covering transformation of $\calE^Q_x$ coincides
with the involution induced by the inversion of the group law.

\item[(iii)] The preimages of $E_{x, 2}$ and $\overline{l}_x$ in $\calE^Q_x$ are irreducible
components of singular fibers. The types of the singular fiber cointainig
the preimages of $E_{x, 2}$ and $\overline{l}_x$  are as follows:

\begin{center}
\begin{tabular}{|c|l|} \hline
$\mbox{I}_2$ & $l_x$ meets $Q$ with two other distinct points \\ \hline
$\mbox{III}$ & $l_x$ is a $3$-fold tangent point. \\ \hline
$\mbox{I}_3$ & $l_x$ is a bitangent line. \\ \hline
$\mbox{IV}$ & $l_x$ is a $4$-fold tangent point. \\ \hline
$\mbox{I}_n$ ($n \ge 4$) & $l_x$ passes through a singular point of type $A_n (n \ge 1)$.  \\ \hline
\end{tabular}
\end{center}
We here use Kodaira's notation (\cite{kodaira}) in order to describe the types of singular fibers.

The following picture describes the case that $l_x$ is  a $3$-fold tangent line at $x$.

%\input pic11.tex 
%WinTpicVersion3.08
\unitlength 0.1in
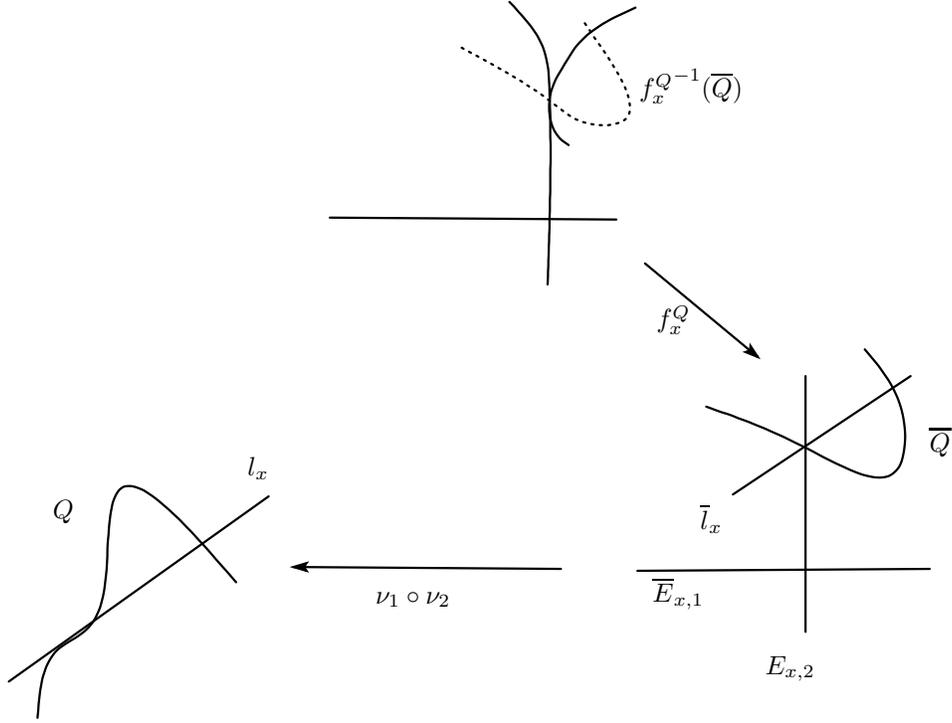
\begin{picture}( 48.2000, 37.5000)(  8.0000,-41.8000)
% LINE 1 0 3 0
% 2 800 3990 2160 3020
% 
\special{pn 13}%
\special{pa 800 3990}%
\special{pa 2160 3020}%
\special{fp}%
% SPLINE 1 0 3 0
% 7 950 4180 1070 3810 1180 3740 1320 3230 1360 3010 1990 3470 1990 3470
% 
\special{pn 13}%
\special{pa 950 4180}%
\special{pa 954 4144}%
\special{pa 956 4108}%
\special{pa 960 4074}%
\special{pa 964 4038}%
\special{pa 970 4006}%
\special{pa 976 3972}%
\special{pa 986 3942}%
\special{pa 996 3912}%
\special{pa 1010 3884}%
\special{pa 1026 3858}%
\special{pa 1044 3834}%
\special{pa 1066 3814}%
\special{pa 1092 3796}%
\special{pa 1120 3778}%
\special{pa 1148 3762}%
\special{pa 1174 3744}%
\special{pa 1198 3726}%
\special{pa 1220 3704}%
\special{pa 1238 3682}%
\special{pa 1254 3656}%
\special{pa 1266 3630}%
\special{pa 1278 3602}%
\special{pa 1286 3572}%
\special{pa 1294 3540}%
\special{pa 1300 3508}%
\special{pa 1304 3474}%
\special{pa 1308 3440}%
\special{pa 1312 3404}%
\special{pa 1314 3368}%
\special{pa 1316 3332}%
\special{pa 1316 3296}%
\special{pa 1318 3258}%
\special{pa 1322 3222}%
\special{pa 1324 3184}%
\special{pa 1328 3148}%
\special{pa 1332 3114}%
\special{pa 1338 3082}%
\special{pa 1344 3052}%
\special{pa 1354 3026}%
\special{pa 1364 3006}%
\special{pa 1376 2988}%
\special{pa 1390 2976}%
\special{pa 1404 2968}%
\special{pa 1422 2966}%
\special{pa 1440 2966}%
\special{pa 1460 2970}%
\special{pa 1482 2978}%
\special{pa 1506 2988}%
\special{pa 1530 3004}%
\special{pa 1554 3020}%
\special{pa 1580 3040}%
\special{pa 1608 3064}%
\special{pa 1636 3088}%
\special{pa 1666 3116}%
\special{pa 1696 3144}%
\special{pa 1726 3176}%
\special{pa 1758 3208}%
\special{pa 1788 3242}%
\special{pa 1822 3278}%
\special{pa 1854 3314}%
\special{pa 1886 3352}%
\special{pa 1920 3390}%
\special{pa 1954 3428}%
\special{pa 1986 3466}%
\special{pa 1990 3470}%
\special{sp}%
% LINE 1 0 3 0
% 2 4090 3410 5620 3400
% 
\special{pn 13}%
\special{pa 4090 3410}%
\special{pa 5620 3400}%
\special{fp}%
% LINE 1 0 3 0
% 2 4970 2390 4970 3730
% 
\special{pn 13}%
\special{pa 4970 2390}%
\special{pa 4970 3730}%
\special{fp}%
% LINE 1 0 3 0
% 2 4590 3010 5520 2390
% 
\special{pn 13}%
\special{pa 4590 3010}%
\special{pa 5520 2390}%
\special{fp}%
% SPLINE 1 0 3 0
% 6 4450 2550 4960 2760 5440 2890 5480 2600 5280 2250 5280 2250
% 
\special{pn 13}%
\special{pa 4450 2550}%
\special{pa 4478 2560}%
\special{pa 4504 2570}%
\special{pa 4532 2580}%
\special{pa 4560 2588}%
\special{pa 4588 2598}%
\special{pa 4616 2610}%
\special{pa 4644 2620}%
\special{pa 4672 2630}%
\special{pa 4700 2642}%
\special{pa 4730 2654}%
\special{pa 4760 2666}%
\special{pa 4790 2680}%
\special{pa 4822 2694}%
\special{pa 4854 2708}%
\special{pa 4886 2724}%
\special{pa 4918 2740}%
\special{pa 4952 2756}%
\special{pa 4988 2774}%
\special{pa 5022 2792}%
\special{pa 5058 2812}%
\special{pa 5094 2830}%
\special{pa 5130 2848}%
\special{pa 5166 2864}%
\special{pa 5200 2880}%
\special{pa 5234 2894}%
\special{pa 5266 2906}%
\special{pa 5298 2914}%
\special{pa 5328 2920}%
\special{pa 5356 2922}%
\special{pa 5382 2920}%
\special{pa 5406 2914}%
\special{pa 5426 2902}%
\special{pa 5444 2886}%
\special{pa 5460 2866}%
\special{pa 5472 2840}%
\special{pa 5480 2810}%
\special{pa 5486 2778}%
\special{pa 5490 2742}%
\special{pa 5492 2706}%
\special{pa 5490 2670}%
\special{pa 5486 2634}%
\special{pa 5480 2598}%
\special{pa 5472 2564}%
\special{pa 5462 2532}%
\special{pa 5450 2500}%
\special{pa 5436 2470}%
\special{pa 5422 2442}%
\special{pa 5406 2414}%
\special{pa 5388 2388}%
\special{pa 5370 2362}%
\special{pa 5352 2336}%
\special{pa 5332 2312}%
\special{pa 5312 2288}%
\special{pa 5292 2264}%
\special{pa 5280 2250}%
\special{sp}%
% LINE 1 0 3 0
% 2 2480 1560 3980 1570
% 
\special{pn 13}%
\special{pa 2480 1560}%
\special{pa 3980 1570}%
\special{fp}%
% SPLINE 1 2 3 0
% 6 3170 670 3660 970 3780 1050 4020 1050 3820 550 3800 520
% 
\special{pn 13}%
\special{pa 3170 670}%
\special{pa 3200 686}%
\special{pa 3228 702}%
\special{pa 3256 716}%
\special{pa 3284 732}%
\special{pa 3314 748}%
\special{pa 3342 764}%
\special{pa 3370 780}%
\special{pa 3398 796}%
\special{pa 3426 812}%
\special{pa 3452 828}%
\special{pa 3480 844}%
\special{pa 3506 862}%
\special{pa 3534 880}%
\special{pa 3560 896}%
\special{pa 3586 914}%
\special{pa 3612 934}%
\special{pa 3636 952}%
\special{pa 3662 972}%
\special{pa 3686 990}%
\special{pa 3712 1010}%
\special{pa 3738 1028}%
\special{pa 3766 1044}%
\special{pa 3798 1058}%
\special{pa 3832 1066}%
\special{pa 3868 1074}%
\special{pa 3904 1076}%
\special{pa 3938 1076}%
\special{pa 3970 1070}%
\special{pa 3998 1062}%
\special{pa 4022 1050}%
\special{pa 4038 1034}%
\special{pa 4048 1014}%
\special{pa 4052 992}%
\special{pa 4050 966}%
\special{pa 4046 938}%
\special{pa 4036 906}%
\special{pa 4024 874}%
\special{pa 4008 842}%
\special{pa 3990 808}%
\special{pa 3968 772}%
\special{pa 3948 738}%
\special{pa 3924 702}%
\special{pa 3902 668}%
\special{pa 3878 636}%
\special{pa 3856 604}%
\special{pa 3836 572}%
\special{pa 3816 544}%
\special{pa 3800 520}%
\special{sp -0.045}%
% VECTOR 1 0 3 0
% 2 3690 3400 2290 3390
% 
\special{pn 13}%
\special{pa 3690 3400}%
\special{pa 2290 3390}%
\special{fp}%
\special{sh 1}%
\special{pa 2290 3390}%
\special{pa 2358 3410}%
\special{pa 2344 3390}%
\special{pa 2358 3370}%
\special{pa 2290 3390}%
\special{fp}%
% STR 2 0 3 0
% 3 1030 3060 1030 3160 2 0
% $Q$
\put(10.3000,-31.6000){\makebox(0,0)[lb]{$Q$}}%
% STR 2 0 3 0
% 3 2050 2830 2050 2930 2 0
% $l_x$
\put(20.5000,-29.3000){\makebox(0,0)[lb]{$l_x$}}%
% STR 2 0 3 0
% 3 2720 3500 2720 3600 2 0
% $\nu_1\circ\nu_2$
\put(27.2000,-36.0000){\makebox(0,0)[lb]{$\nu_1\circ\nu_2$}}%
% STR 2 0 3 0
% 3 4170 3510 4170 3610 2 0
% $\overline{E}_{x,1}$
\put(41.7000,-36.1000){\makebox(0,0)[lb]{$\overline{E}_{x,1}$}}%
% STR 2 0 3 0
% 3 4420 3110 4420 3210 2 0
% $\overline{l}_x$
\put(44.2000,-32.1000){\makebox(0,0)[lb]{$\overline{l}_x$}}%
% STR 2 0 3 0
% 3 5620 2710 5620 2810 2 0
% $\overline {Q}$
\put(56.2000,-28.1000){\makebox(0,0)[lb]{$\overline {Q}$}}%
% STR 2 0 3 0
% 3 4760 3890 4760 3990 2 0
% $E_{x,2}$
\put(47.6000,-39.9000){\makebox(0,0)[lb]{$E_{x,2}$}}%
% STR 2 0 3 0
% 3 4190 2070 4190 2170 2 0
% $f_x^Q$
\put(41.9000,-21.7000){\makebox(0,0)[lb]{$f_x^Q$}}%
% STR 2 0 3 0
% 3 4100 860 4100 960 2 0
% ${f_x^Q}^{-1}(\overline Q)$
\put(41.0000,-9.6000){\makebox(0,0)[lb]{${f_x^Q}^{-1}(\overline Q)$}}%
% VECTOR 1 0 3 0
% 2 4130 1800 4720 2290
% 
\special{pn 13}%
\special{pa 4130 1800}%
\special{pa 4720 2290}%
\special{fp}%
\special{sh 1}%
\special{pa 4720 2290}%
\special{pa 4682 2232}%
\special{pa 4680 2256}%
\special{pa 4656 2264}%
\special{pa 4720 2290}%
\special{fp}%
% SPLINE 1 0 3 0
% 4 3420 430 3600 680 3630 900 3620 1910
% 
\special{pn 13}%
\special{pa 3420 430}%
\special{pa 3444 454}%
\special{pa 3466 480}%
\special{pa 3488 504}%
\special{pa 3508 528}%
\special{pa 3530 554}%
\special{pa 3548 580}%
\special{pa 3566 606}%
\special{pa 3580 634}%
\special{pa 3594 662}%
\special{pa 3604 692}%
\special{pa 3612 722}%
\special{pa 3618 752}%
\special{pa 3624 784}%
\special{pa 3626 816}%
\special{pa 3628 850}%
\special{pa 3630 882}%
\special{pa 3630 914}%
\special{pa 3632 948}%
\special{pa 3632 980}%
\special{pa 3634 1012}%
\special{pa 3634 1046}%
\special{pa 3634 1078}%
\special{pa 3634 1110}%
\special{pa 3634 1142}%
\special{pa 3636 1174}%
\special{pa 3636 1206}%
\special{pa 3636 1238}%
\special{pa 3636 1270}%
\special{pa 3634 1302}%
\special{pa 3634 1336}%
\special{pa 3634 1366}%
\special{pa 3634 1398}%
\special{pa 3634 1430}%
\special{pa 3632 1462}%
\special{pa 3632 1494}%
\special{pa 3632 1526}%
\special{pa 3630 1558}%
\special{pa 3630 1590}%
\special{pa 3630 1622}%
\special{pa 3628 1654}%
\special{pa 3628 1684}%
\special{pa 3626 1716}%
\special{pa 3626 1748}%
\special{pa 3624 1780}%
\special{pa 3624 1812}%
\special{pa 3622 1844}%
\special{pa 3622 1874}%
\special{pa 3620 1906}%
\special{pa 3620 1910}%
\special{sp}%
% SPLINE 1 0 3 0
% 6 4080 460 3780 650 3700 770 3640 890 3700 1160 3730 1180
% 
\special{pn 13}%
\special{pa 4080 460}%
\special{pa 4050 474}%
\special{pa 4020 488}%
\special{pa 3990 502}%
\special{pa 3960 516}%
\special{pa 3932 532}%
\special{pa 3904 548}%
\special{pa 3878 566}%
\special{pa 3852 584}%
\special{pa 3826 604}%
\special{pa 3804 626}%
\special{pa 3782 648}%
\special{pa 3762 672}%
\special{pa 3744 700}%
\special{pa 3728 726}%
\special{pa 3710 754}%
\special{pa 3694 782}%
\special{pa 3678 808}%
\special{pa 3662 836}%
\special{pa 3650 866}%
\special{pa 3638 898}%
\special{pa 3632 930}%
\special{pa 3628 964}%
\special{pa 3628 1000}%
\special{pa 3632 1034}%
\special{pa 3640 1068}%
\special{pa 3652 1100}%
\special{pa 3670 1128}%
\special{pa 3690 1152}%
\special{pa 3716 1170}%
\special{pa 3730 1180}%
\special{sp}%
\end{picture}%

%%%%%%%%%%%%%%%%%%%%%%
%
%'±'±'ɐ}'ð'¢'ê'éD
%
%%%%%%%%%%%%%%%%%%%%

\item[(iv)] Other singular fibers of $\calE_x^Q$ correspond to lines in $\Lambda_x$ not
meeting $Q$ at $4$ distinct points. We refer to \cite[Table 6.2]{miranda-persson} for details. 

\end{enumerate}

\begin{rem}\label{rem:quartic-elli}{\rm Note that  any rational elliptic surface ${\mathcal E}$ with 
at least one reducible singular fiber is obtained above. Namely ${\mathcal E} = {\mathcal E}_x^Q$ for
some $Q$ and a smooth point $x$ on $Q$. 
}
\end{rem}
\subsection{The Mordell-Weil lattices of $\calE_x^Q$}

In this subsection, we give a table of types of singularities of $Q$, the relative position of 
$l_x$ and $Q$, and the Mordell-Weil lattices of $\calE_x^Q$.  We first note that
$\MW(\calE_x^Q)$ has no $2$-torsion, since we assume that $Q$ is irreducible. 
Also we
omit cases which never occur. As for the  structure of the Mordell-Weil lattices for rational elliptic surfaces, we refer to 
\cite{oguiso-shioda} and  
to \cite{shioda92}  for the correction of the misprints in \cite{oguiso-shioda}. Let us explain notations used in the table.

%\begin{itemize}

$\bullet$ $\Xi_Q$ and  $l_x\cap Q$ are the same as those in the table Theorem~\ref{thm:main}

$\bullet$ $R_{Q,x}$: the subgroup of $\NS(\calE^Q_x)$ generated by 
$\{\Theta_{v,1}, \ldots, \Theta_{v,m_v-1}\}_{v \in R_{\varphi_x^Q}}$. Note tha $R_{Q,x}$ is 
isomorphic to a direct sum of root lattices of A-D-E type, and we describe $R_{Q, x}$
as a direct sum of them.  

$\bullet$ $\MW(\calE_x^Q)$: the lattice strucure of $\MW(\calE_x^Q)$. To describe them, we
use the notation in \cite{oguiso-shioda}. Namely $\bullet^*$ means the dual lattice of the
lattice $\bullet^*$ and $\langle m \rangle$ denotes a lattice of rank $1$, $\ZZ x$ with
$\langle x, x \rangle = m$. Also a matrix means the intersection matrix with respect to
a certain basis. Note that the lattice structure is determined by
$R_{Q,x}$ as $\MW(\calE_x^Q)$ has no $2$-torsion. 

$\bullet$ $\MW^0(\calE_x^Q)$: the narrow part of $\MW(\calE_x^Q)$, i.e., the subgroup of
$\MW(\calE_x^Q)$ consisting of  sections $s$ with $s\Theta_{v,0} = 1$.

%\end{itemize}

\begin{center}
\begin{tabular}{|c|c|c|c|c|c|} \hline 
 No. & $\Xi_Q$ & $l_x\cap Q$ & $R_{Q,x}$ & $\MW(\calE_x^Q)$ & $\MW^0(\calE_x^Q)$ \\ \hline
 1  & $A_6$    &     $s$   & $A_6\oplus A_1$ & $\langle 1/14 \rangle$ & $\langle 14 \rangle$ \\ 
 2 & $A_6$    &      $sb$  & $A_8$                 & $\ZZ/3\ZZ$ & $\{0\}$ \\ \hline
 3 & $E_6$    &      $s$    & $E_6\oplus A_1$ & $\langle 1/6 \rangle$ & $\langle 6 \rangle$ \\
 4 & $E_6$    &      $b$   & $E_6\oplus A_2$ & $\ZZ/3\ZZ$ & $\{0\}$  \\ \hline    
 5 & $A_5$   &       $s$   & $A_5\oplus A_1$ & $A_1^*\oplus \langle 1/6 \rangle$ & $A_1\oplus 
 \langle 6 \rangle$ \\
 6 & $A_5$    &      $b$   & $A_5\oplus A_2$ & $A_1^* \oplus \ZZ/3\ZZ$ & $A_1$ \\
 7 & $A_5$    &      $sb$  & $A_7$                & $\langle 1/8 \rangle$   &$\langle 8 \rangle$ \\ \hline
  8 & $D_5$    &      $s$    & $D_5\oplus A_1$ & $A_1^*\oplus \langle 1/4\rangle$ & $A_1\oplus \langle 4 \rangle$ \\
   9 & $D_5$    &      $b$    & $D_5\oplus A_2$ & $\langle 1/12 \rangle$  & $\langle 12 \rangle$ \\ \hline
 10 & $D_4$    &      $s$    & $D_4\oplus A_1$ & $(A_1^*)^{\oplus 3}$ & $A_1^{\oplus 3}$ \\
 11 & $D_4$    &      $b$    & $D_4\oplus A_2$ & $\displaystyle 
                                                                       \frac 16 \left (
                                                                          \begin{array}{cc}
                                                                          2 & 1 \\
                                                                          1 & 2  
                                                                          \end{array} \right )$ & 
                                                                          $\displaystyle \left (
                                                                          \begin{array}{cc}
                                                                          4 & -2 \\
                                                                          -2 & 4  
                                                                          \end{array} \right )$ \\ \hline
 12 & $A_4 + A_2$ &$s$ &$A_4\oplus A_2\oplus A_1$ & $\langle 1/30 \rangle$ & $\langle 30 \rangle$ \\
  13 & $A_4 + A_2$ &$sb$ &$A_4\oplus A_4$ & $\ZZ/5\ZZ$ & $\{0 \}$ \\ \hline
  14 & $A_4 + A_1$& $s$ & $A_4\oplus A_1^{\oplus 2}$ & $\displaystyle 
                                                                       \frac 1{10} \left (
                                                                          \begin{array}{cc}
                                                                          2 & 1 \\
                                                                          1 & 3  
                                                                          \end{array} \right )$ & 
                                                                          $\displaystyle \left (
                                                                          \begin{array}{cc}
                                                                          6 & -2 \\
                                                                          -2 & 4  
                                                                          \end{array} \right )$ \\
  15 & $A_4 + A_1$ &   $b$ & $A_4\oplus A_2 \oplus A_1$ & $\langle 1/30 \rangle$ &
    $\langle 30 \rangle$ \\
16 & $A_4 + A_1$ & $sb$ & $A_4 \oplus A_3$ & $\langle 1/20 \rangle$ & $\langle 20 \rangle$ \\
17 & $A_4 + A_1$&  $sb$ & $A_6 \oplus A_1$ & $\langle 1/14 \rangle$ & $\langle 14 \rangle$ \\ \hline    
18 & $A_3 + A_2$ & $s$ &$A_3\oplus A_2 \oplus A_1$ & $A_1^* \oplus \langle 1/12 \rangle$  &
$A_1 \oplus \langle 12 \rangle$ \\   
19 & $A_3 + A_2$ & $sb$ & $A_4 \oplus A_3$ & $\langle 1/20 \rangle$ & $\langle 20 \rangle$  \\
20 & $A_3 + A_2$ & $sb$ & $A_5 \oplus A_2$ & $A_1^*\oplus \ZZ/3\ZZ$ & $A_1$ \\ \hline 
21 & $A_3 + A_1$ & $s$ & $A_3 \oplus A_1^{\oplus 2}$ & $(A_1^*)^{\oplus 2}\oplus \langle 1/4 \rangle$ & $A_1^{\oplus 2}\oplus \langle 4 \rangle$ \\                                                              
22 & $A_3 + A_1$ & $b$ & $A_3\oplus A_2 \oplus A_1$ & $A_1^*\oplus \langle 1/12 \rangle$ &
$A_1\oplus \langle 12 \rangle$ \\   
23 & $A_3 + A_1$ & $sb$ & $A_5 \oplus A_1$ & $A_1^*\oplus \langle 1/6 \rangle$ &
$A_1\oplus \langle 6 \rangle$ \\   
24  & $A_3 + A_1$ & $sb$ & $A_3 \oplus A_3$ & $ \langle 1/4 \rangle^{\oplus 2}$ &
$ \langle 4 \rangle^{\oplus 2}$ \\   \hline
25 & $3A_2$ & $s$ & $A_2^{\oplus 3}\oplus A_1$ & $\langle 1/6 \rangle\oplus \ZZ/3\ZZ$ &
$\langle 6 \rangle$ \\
26 & $3A_2$ & $b$ & $A_2^{\oplus 4}$ & $(\ZZ/3\ZZ)^{\oplus 2}$ & $\{0\}$ \\ \hline
27 & $2A_2 + A_1$ & $s$ & $A_2^{\oplus 2}\oplus A_1^{\oplus 2}$ & $\langle 1/6 \rangle^{\oplus 2}$ & $\langle 6 \rangle^{\oplus 2}$ \\
28 & $2A_2 + A_1$ & $b$ & $A_2^{\oplus 3}\oplus A_1$ & $\langle 1/6\rangle\oplus \ZZ/3\ZZ$ &
$\langle 6 \rangle$ \\
29 & $2A_2 + A_1$ & $sb$ & $A_4\oplus A_2 \oplus A_1$ & $\langle 1/30 \rangle$ &
$\langle 30 \rangle$ \\ \hline
30 & $A_2 + 2A_1$ & $s$ & $A_2\oplus A_1^{\oplus 3}$ & $A_1^*\oplus  \displaystyle 
                                                                       \frac 16 \left (
                                                                          \begin{array}{cc}
                                                                          2 & 1 \\
                                                                          1 & 2  
                                                                          \end{array} \right )$ & 
                                                                          $A_1\oplus \displaystyle \left (
                                                                          \begin{array}{cc}
                                                                          4 & -2 \\
                                                                          -2 & 4  
                                                                          \end{array} \right )$ \\ 
 31& $A_2 + 2A_1$ & $b$ & $A_2^{\oplus 2}\oplus A_1^{\oplus 2}$ & $\langle 1/6 \rangle^{\oplus 2}$  & $\langle 6 \rangle^{\oplus 2}$ \\
 32 & $A_2 + 2A_1$ & $sb$ & $A_4 \oplus A_1^{\oplus 2}$ & $\displaystyle 
                                                                       \frac 1{10} \left (
                                                                          \begin{array}{cc}
                                                                          2 & 1 \\
                                                                          1 & 3  
                                                                          \end{array} \right )$ & 
                                                                          $\displaystyle \left (
                                                                          \begin{array}{cc}
                                                                          6 & -2 \\
                                                                          -2 & 4  
                                                                          \end{array} \right )$ \\ 
 33 & $A_2 + 2A_1$ & $sb$ & $A_3\oplus A_2 \oplus A_1$ & $A_1^*\oplus \langle 1/12 \rangle$ &
 $A_1\oplus \langle 12 \rangle$ \\ \hline      

 34 & $3A_1$ & $s$ & $A_1^{\oplus 4}$ & $(A_1^*)^{\oplus 4}$ & $A_1^{\oplus 4}$ \\
 35 & $3A_1$ & $s$ & $A_2\oplus A_1^{\oplus 3}$ &   $A_1^*\oplus  \displaystyle 
                                                                       \frac 16 \left (
                                                                          \begin{array}{cc}
                                                                          2 & 1 \\
                                                                          1 & 2  
                                                                          \end{array} 
                                                                          \right )$ & 
                                                                          $A_1\oplus \displaystyle \left (
                                                                          \begin{array}{cc}
                                                                          4 & -2 \\
                                                                          -2 & 4  
                                                                          \end{array} \right )$ \\                                                                                                                                                                                                                
36 & $3A_1$ & $sb$ & $A_3\oplus A_1^{\oplus  2}$ & $(A_1^*)^{\oplus 2}\oplus \langle 1/4\rangle $ &
$A_1^{\oplus 2}\oplus \langle 4 \rangle$ \\ \hline  
 \end{tabular}
 
 \begin{tabular}{|c|c|c|c|c|c|} \hline 
 No. & $\Xi_Q$ & $l_x\cap Q$ & $R_{Q,x}$ & $\MW(\calE_x^Q)$ & 
 $\MW^0(\calE_x^Q)$ \\ \hline
                                                             
 37 & $A_4$ & $s$ &$A_4\oplus A_1$ & $\frac 1{10} \left (
                                                                          \begin{array}{ccc}
                                                                          3 & 1 & -1  \\
                                                                          1 & 7 & 3 \\
                                                                          -1 & 3 & 7  
                                                                          \end{array} 
                                                                          \right )$ & 
                                                                          $\displaystyle \left (
                                                                          \begin{array}{ccc}
                                                                          4 & -1  & 1  \\
                                                                          -1 & 2 & -1 \\
                                                                          1 & -1 & 2  
                                                                          \end{array} \right )$ \\ 
 38 & $A_4$ & $b$ &$A_4\oplus A_2$ &  $\frac 1{15} \left (
                                                                          \begin{array}{cc}
                                                                          2 & 1   \\
                                                                          1 & 8 \\
                                                                           \end{array} 
                                                                          \right )$ & 
                                                                          $\displaystyle \left (
                                                                          \begin{array}{cc}
                                                                          8 & -1  \\
                                                                          -1 & 2  
                                                                           \end{array} \right )$ \\ 
                                                            
39 &  $A_4$ & $sb$ & $A_6$                   &  $\frac 1{7} \left (
                                                                          \begin{array}{cc}
                                                                          2 & 1   \\
                                                                          1 & 4 \\
                                                                           \end{array} 
                                                                          \right )$ & 
                                                                          $\displaystyle \left (
                                                                          \begin{array}{cc}
                                                                          4 & -1  \\
                                                                          -1 & 2  
                                                                           \end{array} \right )$ \\ \hline
40 & $A_3$ &  $s$ & $A_3\oplus A_1$ & $A_3^*\oplus A_1^*$ & $A_3\oplus A_1$ \\   
41 & $A_3$ & $b$ & $A_3\oplus A_2$ &  $\frac 1{12} \left (
                                                                          \begin{array}{ccc}
                                                                          7 & 1 & 2  \\
                                                                          1 & 7 & 2 \\
                                                                          2 & 2 & 4  
                                                                          \end{array} 
                                                                          \right )$ & 
                                                                          $\displaystyle \left (
                                                                          \begin{array}{ccc}
                                                                          2 & 0  & -1  \\
                                                                          0 & 2 & -1 \\
                                                                          -1 & -1 & 4  
                                                                          \end{array} \right )$ \\ 
42 & $A_3$ & $sb$ & $A_5$ & $A_2^*\oplus A_1^*$ & $A_2\oplus A_1$ \\ \hline
43 & $2A_2$ & $s$ & $A_2^{\oplus 2}\oplus A_1$ & $A_2^*\oplus \langle 1/6 \rangle$ &
$A_2\oplus \langle 6 \rangle$ \\
44 & $2A_2$ & $b$ & $A_2^{\oplus 3}$ & $A_2^*\oplus \ZZ/3\ZZ$ & $A_2$ \\
45 & $2A_2$ & $sb$ & $A_4\oplus A_2$ &    $\frac 1{15} \left (
                                                                          \begin{array}{cc}
                                                                          2 & 1   \\
                                                                          1 & 8 \\
                                                                           \end{array} 
                                                                          \right )$ & 
                                                                          $\displaystyle \left (
                                                                          \begin{array}{cc}
                                                                          8 & -1  \\
                                                                          -1 & 2  
                                                                           \end{array} \right )$ \\ \hline
                                                                        
46 & $A_2 + A_1$ & $s$ & $A_2\oplus A_1^{\oplus 2}$ & $\frac 1{6} \left (
                                                                          \begin{array}{cccc}
                                                                          2 & 1 & 0 & -1  \\
                                                                          1 & 5 & 3 & 1 \\
                                                                          0 & 3 & 6 & 3 \\
                                                                          -1 & 1 & 3 & 5  
                                                                          \end{array} 
                                                                          \right )$ & 
                                                                          $\displaystyle \left (
                                                                          \begin{array}{cccc}
                                                                          4 & -1  & 0 & 1  \\
                                                                          -1 & 2 & -1 & 0 \\
                                                                          0 & -1 & 2  & -1 \\
                                                                          1 & 0 & -1 & 2  
                                                                          \end{array} \right )$ \\ 
 47 & $A_2 + A_1$ & $b$ & $A_2^{\oplus 2} \oplus A_1$ & $A_2^*\oplus \langle 1/6\rangle$ &
 $A_2\oplus \langle 6 \rangle$ \\
 48 & $A_2 + A_1$ & $sb$ & $A_4\oplus A_1$ & $\frac 1{10} \left (
                                                                          \begin{array}{ccc}
                                                                          3 & 1 & -1  \\
                                                                          1 & 7 & 3 \\
                                                                          -1 & 3 & 7  
                                                                          \end{array} 
                                                                          \right )$ & 
                                                                          $\displaystyle \left (
                                                                          \begin{array}{ccc}
                                                                          4 & -1  & 1  \\
                                                                          -1 & 2 & -1 \\
                                                                          1 & -1 & 2  
                                                                          \end{array} \right )$ \\ 
  49 &  $A_2 + A_1$ & $sb$ & $A_4\oplus A_1$ & $\frac 1{12} \left (
                                                                          \begin{array}{ccc}
                                                                          7 & 1 & 2  \\
                                                                          1 & 7 & 2 \\
                                                                          2 & 2 & 4  
                                                                          \end{array} 
                                                                          \right )$ & 
                                                                          $\displaystyle \left (
                                                                          \begin{array}{ccc}
                                                                          2 & 0  & -1  \\
                                                                          0 & 2 & -1 \\
                                                                          -1 & -1 & 4  
                                                                          \end{array} \right )$ \\ \hline
                                                                          
 50 & $2A_1$ & $s$ & $A_1^{\oplus 3}$ & $D_4^* \oplus A_1^*$ & $D_4\oplus A_1$ \\                                                                       
51 & $2A_1$ & $b$  & $A_2\oplus A_1^{\oplus 2}$ & $\frac 1{6} \left (
                                                                          \begin{array}{cccc}
                                                                          2 & 1 & 0 & -1  \\
                                                                          1 & 5 & 3 & 1 \\
                                                                          0 & 3 & 6 & 3 \\
                                                                          -1 & 1 & 3 & 5  
                                                                          \end{array} 
                                                                          \right )$ & 
                                                                          $\displaystyle \left (
                                                                          \begin{array}{cccc}
                                                                          4 & -1  & 0 & 1  \\
                                                                          -1 & 2 & -1 & 0 \\
                                                                          0 & -1 & 2  & -1 \\
                                                                          1 & 0 & -1 & 2  
                                                                          \end{array} \right )$ \\ 
52 & $2A_1$ & $sb$ & $A_3\oplus A_1$ & $A_3^* \oplus A_1^*$ & $A_3\oplus A_1$ \\ \hline  
53 & $A_2$ & $s$ & $A_2\oplus A_1$ & $A_5^*$ & $A_5$ \\
54 & $A_2$ & $b$ & $A_2^{\oplus 2}$ & $(A_2^*)^{\oplus 2}$ & $A_2^{\oplus 2}$ \\
55 & $A_2$ & $sb$ & $A_4$ & $A_4^*$ & $A_4$ \\ \hline

%\end{tabular}

%
%\begin{tabular}{|c|c|c|c|c|c|}\hline
%No. & $\Xi_Q$ & $l_x\cap Q$ & $R_{Q,x}$ & $\MW$ & $\MW^0$ \\ \hline
56 & $A_1$ & $s$ & $A_1^{\oplus 2}$ & $D_6^*$ &$D_6$ \\
57 & $A_1$ & $b$ & $A_2\oplus A_1$ & $A_5^*$ & $A_5$ \\
58 & $A_1$ & $sb$ & $A_3$ & $D_5^*$ & $D_5$ \\ \hline   
59 & $\emptyset$ & $s$ &$A_1$ & $E_7^*$ & $E_7$ \\
60 & $\emptyset$ & $b$ & $A_2$ & $E_6^*$ & $E_6$ \\ \hline                                                                                                                                       
\end{tabular}                                                                                             
\end{center}

%
%%
%
%
%

%%%%%%%%%%%%%%%%%%%%%%%%%%%%%%%%%%%%% 
%  \input sec-coqua2.tex
 %%%%%%%%%%%%%%%%%%%%%%%%%%%%%%%%%%%%%% 
 %%%%%%%%%%%%%%%%%%%%%%%%%%%%%%%%%%%%%%%

\section{Proof of Theorem~\ref{thm:main}} 

We keep the same notations as before.
Our proof is done by a case-by-case checking. Let us start with the following lemma.

\begin{lem}\label{lem:etc-sec}{Let $C$ be an even tangential conic of $Q$ through $x$.
The preimage of $C$ in ${\calE}_x^Q$ consists of two sections $s_C^+$ and $s_C^-$ such that 

\medskip

$(i)$ $\langle s_C^+, s_C^+ \rangle = \langle s_C^-, s_C^-\rangle = 2$

$(ii)$ $s_C^+O = s_C^-O = 0$

$(iii)4$ $s_C^+\Theta_{v,0} = s_C^-\Theta_{v,0} = 1$ {\it for all} $v \in R_{\varphi_x^Q}$, i.e,
$s_C^{\pm} \in \MW^0(\calE_x^Q)$.

Coversely, for any section $s$ in $\MW(\calE_x^Q)$ with two of the above three properties,
the image of $s$ in $\PP^2$ is an even tangential conic to $Q$.
}
\end{lem}

\proof 
We first note that two of the properties $(i)$, $(ii)$ and $(iii)$ imply the remaining. This follows from
a formula
\[
\langle s, s \rangle = 2 + 2sO  - \sum_{v \in R_{\varphi}}
\mbox{Corr}_v(s, s).
\]
for a rational elliptic surface $\calE$ and $s \in \MW(\calE_x^Q)$.

Let $\overline{C}$ be the proper transform of $C$ in $\widehat{\PP}^2$. Since
$\overline{C}$ is tangent to $\overline{Q}$ at each intersection point and 
$\overline{C}\cap \overline{E}_{x,1} = \emptyset$, the preimage of $\overline {C}$ in
$\calE^Q_x$ consists of $2$ irreducible components $s_C^+$ and $s_C^-$ so that
$s_C^{\pm}O = 0$. Since $\overline{C}$ meets the propertransform of a general 
member in
$\Lambda_x$ at one point, both $s_C^+$ and $s_C^-$ are sections of 
$\varphi_x^Q : \calE^Q_x \to \PP^1$.  The property $(iii)$ follows from the fact that
$\overline{C}$ meets $E_{x,2}$ and $\overline{C}$ does not pass through 
singularities of $\overline{Q}$. Now the property $(i)$ is straightforward from the
explicit formula for $\langle\, , \, \rangle$.

Conversely, suppose that we have a section $s$  with two of the properties $(i)$, $(ii)$ and $(iii)$. Let  $C_s$ be the image
of $s$ in $\PP^2$. By our construction of ${\calE}_x^Q$, we infer that $C_s$ is a conic tangent to $Q$ at
$x$.  Since $C_s$ is also the image of $\sigma_{f_x^Q}^*s$, we infer that $C_s$ is an even tangent conic to $Q$.
\proofend

\begin{thm}\label{thm:qrc}{ Let $C$ be an even tangential conic to $Q$ and let $s_C^+$ be the section as above.
\[
(C/Q) = (-1)^{\varepsilon(s_C^+)},
\]
where the symbol $\varepsilon(s_C^+)$ is the same as that defined in Theorem~\ref{thm:recipro}.
}
\end{thm}

\proof Let $\widehat{\PP}^2$ as before. Since $\overline{l}_x$ is a $(-1)$ curve, by blowing down 
$\overline{l}_x$, we obtain $\Sigma_2$ with the following properties:

\begin{enumerate}

\item[(i)] The image of $\overline{Q}$ is a trisection $\Gamma_Q \sim 3(\Delta_0 + 2F)$.

\item[(ii)] Singularities of $\Gamma_Q$ is the same as those of $Q$ except the $A_1$ singularity
 caused by blowing down $\overline{l}_x$.

\item[(iii)] The image of $\overline{E}_{x,1} = \Delta_0$.

\item[(iv)] The image of $\overline{C}$ is a section $\Delta_C$ such that $\Delta_C \sim
(\Delta_0 + 2F)$ and $\Delta_C$ is even tangent to $\Gamma_Q$. 

\end{enumerate}

Let $f_o : Z_o \to \Sigma_2$  be the induced double cover by $f_C : Z_C \to \PP^2$, i.e.,
the $\CC(Z_C)$-normalization  of $\Sigma_2$. One easily see that $\Delta_{f_o} = \Delta + \Delta_C$.

\[
\begin{diagram}
\node{Z_C}\arrow{s} \node{\calE_x^Q}\arrow{s} \node{Z_o}\arrow{s}\\
\node[1]{\PP^2} \node[1]{\widehat{\PP}^2}\arrow{w}\arrow{e}\node[1]{\Sigma_2}
\end{diagram}
\]

Since $\Delta_C$ is the image of $\overline {C}$, it is also the image of $s_C^{\pm}$. Hence we infer that
\[
(C/Q) = 1 \Leftrightarrow (\Delta_0 + \Delta_C/\Gamma_Q) = 1.
\]

Hence by Theorem~\ref{thm:recipro}, we infer that $(C/Q) = 1$ if and only if
$s_C^+ \in 2\MW(\calE_x^Q)$.  \proofend

\begin{rem}\label{cor:qrc}{\rm For $s_o$ in Theorem~\ref{thm:qrc}, we have $\langle s_o, s_o \rangle = 1/2$.  Hence
if $\MW(\calE_x^Q)$ has no section $s$ with $\langle s, s\rangle = 1/2$, there is no quadratic
residue even tangential conic through $x$.
}
\end{rem}

\begin{lem}\label{lem:genus}{ Let $\widetilde Q$ be the normalization of $Q$ and we denote the genus of 
$\widetilde Q$ by $g(\widetilde Q)$.

$(i)$ There exists no even tangential conic to $Q$ which is quadratic residue $\bmod Q$ if 
$g(\widetilde {Q}) \ge 2$.

$(ii)$ All even tangential conic to $Q$ are quadratic residue $\bmod Q$ if $g(\widetilde {Q}) = 0$. 
}
\end{lem}

\proof $(i)$ Let $C$ be an even tangential conic to $Q$ and suppose that $(C/Q) = 1$.  Let $f_C : Z_C \to \PP^2$ be a double cover with $\Delta_{f_C} = C$. Then 
$f_C^*Q$ is of the form $Q^+ + Q^-$. Since $\PP^1\times\PP^1, \Pic(Z_C)\cong \ZZ\oplus \ZZ$ 
and the covering transformation induces an involution $(a, b) \mapsto (b, a)$ on $\Pic(Z_C)$, we
infer that $Q^+ \sim Q^- \sim (2, 2)$. Since $Q^+, Q^-$ and $Q$ are birationally equivalent, we
have $g(\widetilde {Q}) \le 1$. This leads us to a contradiction.

$(ii)$ Since the induced double cover on $\widetilde{Q}$ is unramified, $(C/Q) = 1$.
\proofend

 Now we easily have the following theorem:

\begin{thm}\label{thm:key}{ Let $Q$ be an irreducible quartic. Choose a smooth point $x \in Q$ and
let $\calE_x^Q$ be the rational elliptic surface as in \S 1. Then we have the following:

$(i)$ Let $\mbox{\rm ETC}$ be the set of conics passing through $x$. Then
\begin{eqnarray*}
 \sharp \mbox{\rm ETC} &=& \sharp\{s \in \MW(\calE_x^Q) \mid \langle s, s \rangle = 2, sO = 0\}/2 \\
 & = & \sharp\{s \in \MW^0(\calE_x^Q) \mid \langle s, s \rangle = 2\}/2
\end{eqnarray*}

$(ii)$  Let $\mbox{\rm QRETC}$ be the set of even tangential conics passing through $x$ with $(C/Q) = 1$. Then
\begin{eqnarray*}
\sharp \mbox{\rm QRETC} & =  & \sharp\{s \in 2\MW(\calE_x^Q) \mid  \langle s, s \rangle = 2, sO = 0\}/2 \\
& = & \sharp\{ s \in 2\MW(\calE_x^Q)\cap \MW^0(\calE_x^Q) \mid \langle s, s \rangle = 2 \}/2
\end{eqnarray*}
}
\end{thm}

\proof  Our statements $(i)$ and $(ii)$ are immediate from Lemma~\ref{lem:etc-sec} and
Theorem~\ref{thm:qrc}.  
\proofend

We now prove Theorem~\ref{thm:main}  by a case-by-case checking. We first compute
$\sharp\mbox{\rm ETC}$.
By Lemma~\ref{lem:etc-sec}, it is enough to see the number of sections $s$ in 
the narrow part $\MW^0(\calE_x^Q)$ of $\MW(\calE_x^Q)$ with $\langle s, s \rangle =2$. 

For the lattices of A-D-E types, it is nothing but the number of roots, and the following table
is well-known (see \cite{conway-sloane})

\begin{center}
\begin{tabular}{cccc} \hline
$A_n$ & $D_n$ ($n \ge 4$) & $E_6$ & $E_7$ \\ \hline
$n(n+1)$ & $2n(n-1)$ &  $72$ & $126$ \\ \hline
\end{tabular}
\end{center}

From the above table and that in \S 2, our statement on $\sharp \mbox{ETC}$ is straightforward except
for the cases $11, 14, 30, 32, 35, 37, 38, 39, 41, 45, 46, 48, 49, 51$.  For the rank $2$ cases
among the exceptional
cases, our statement follows easily by direct computation. For the cases of $\rank  > 2$, we make
use of \cite[Lemma 3.8]{oguiso-shioda}, which is as follows:
\[
\left ( \begin{array}{ccc}
          4 & -1 & 1 \\
          -1 & 2 & -1 \\
          1 & -1 & 2 \end{array} \right ) \cong \mbox{$A_1^{\perp}$ in $A_4$}, \quad 
 \left ( \begin{array}{cccc}
          4 & -1 & 0& 1 \\
          -1 & 2 & -1 & 0 \\
          0 & -1 & 2 & -1 \\
          1 &0 &  -1 & 2 \end{array} \right ) \cong \mbox{$A_1^{\perp}$ in $A_5$}
\]
\[
\left ( \begin{array}{ccc}
          2 & 0 & -1 \\
          0 & 2 & -1 \\
          -1 & -1 & 2 \end{array} \right ) \cong \mbox{$A_2^{\perp}$ in $D_5$}
\]
where the terminology $\bullet^{\perp}$ in $\blacksquare$ means that we embed a lattice $\bullet$ into
$\blacksquare$ and $\bullet^{\perp}$ is the orthogonal complement of $\bullet$ in 
$\blacksquare$.  Also, by \cite[Lemma 3.8]{oguiso-shioda}, the embedding is determined up to isomoprhisms.  Hence we just count the number of 
roots which are orthogonal to the embedded lattices. To be more precise, we explain the case
$A_1^{\perp}$ in $A_5$. We first consider the realization of $A_5$ as follows:
\[
A_5 = \{(x_1, \ldots, x_6) \mid \sum_i x_i= 0, x_i \in \ZZ\} \subset \RR^6
\]
and the pairing is induced from the Euclidean metric $\sum_i x_i^2$ in $\RR^6$. Under these
circumstances,
the roots are given by a vector $(1, -1, 0,0,0,0)$ and those obtained by permutations of the 
coordinates. We fix an embedding of $A_1$ given by $\ZZ(1, -1, 0,0,0,0) \subset A_5$. Then roots in
$A_1^{\perp}$ are 
\[
\begin{array}{ccc}
 (0,0, \pm 1, \mp 1, 0,0) & (0,0, \pm 1, 0, \mp 1, 0) &
 (0, 0, \pm 1, 0, 0, \mp 1) \\
   (0,0,0, \pm 1, \mp 1, 0) & (0,0,0, \pm 1, 0, \mp 1) & (0,0,0,0, \pm 1, \pm 1) . 
\end{array}
\]
 Since the remaining cases are similar, we omit them. Thus we have a list for $\sharp 
 \mbox{\rm ETC}$.
 
 We now go on to compute $\sharp \mbox{\rm QRETC}$. 
 We first note that $\sharp\mbox{\rm QRETC} = 0$ if
 $\sharp\mbox{\rm ETC} = 0$. In the following, we only cosider the case of $\sharp\mbox{\rm ETC} \neq 0$. 
 
 Since $Q$ is irreducible, $\MW(\calE_x^Q)$ has
 no $2$-torsion. Hence for each $s \in 2\MW(\calE_x^Q)$, there exists a unique $s_o \in \MW(\calE_x^Q)$
 such that $2s_o = s$.  For distinct $C_1, C_2 \in \mbox{\rm QRETC}$, $s_{C_1}^+$ and $s_{C_2}^+$ are
 distinct in $\MW^0(\calE_x^Q)$. Hence it is enough to compute
 \[
 \sharp \{ s_o \in \MW(\calE_x^Q) \mid \langle s_o, s_o\rangle = 1/2, 2s_o \in \MW^0(\calE_x^Q) \}
 \]
 Now Theorem~\ref{thm:main} follows from the following claim:

 \medskip
 
 \textbf{Claim.} {\it Suppse that $\sharp\mbox{\rm ETC} \neq 0$. If $\MW(\calE_x^Q)$ has an $A_1^*$ as a direct summand. then two generators
 $\pm {\tilde s}$ of $A_1^*$ are sections such that $\langle \tilde s, \tilde s\rangle = 1/2, 
 2\tilde s \in \MW^0(\calE_x^Q)$.
 Conversely if there exists $s_o \in \MW(\calE_x^Q)$ such that $\langle s_o, s_o\rangle = 1/2, 2s_o \in \MW^0(\calE_x^Q)$,
 then  $\ZZ s_o (\cong A_1^*)$ is a direct summand of $\MW(\calE_x^Q)$.
 }

 \medskip

  \textsl{Proof of Claim.} Suppose that $A_1^*$ is a direct summand of $\MW(\calE_x^Q)$ and
  let $\tilde s$ be a section such that $\ZZ{\tilde s} = A_1^*$.  Then $\langle \tilde s, \tilde s \rangle
  = 1/2$ and $2\tilde s \in \MW^0(\calE_x^Q)$ by \cite[Theorem 9.1]{shioda90}.

We now go on to show the converse.   Let $s_o$ be a section with $\langle s_o, s_o \rangle = 1/2,
2s_o \in \MW^0(\calE_x^Q)$. 
 As for the dual lattices of A-D-E type, we have the following table:
 
 \begin{center}
\begin{tabular}{c|cccc} \hline
 Type &$A_n^*$ &  $D_n^*$ ($n \ge 4$) & $E_6^*$ & $E_7^*$ \\ \hline
Minimum norm& $\frac{n}{(n+1)}$ & 1  &  $\frac 43 $ & $\frac 32 $ \\ \hline
\end{tabular}
\end{center}

%Hence $\MW(\calE_x^Q)$ has an $A_1^*$ as a direct summand except for the cases 11, 14, 30, 32, 35, 37, 38, 39, 41, 45, 46, 48,
%49,51.

 Hence we easily see that $\MW(\calE_x^Q)$ has an $A_1^*$ direct summand except for the cases
37, 38, 39, 41, 45, 46, 48,
49,51. We see that there is no section $s$ with $\langle s, s \rangle = 1/2$ for these
exceptional cases.

\medskip
\textsl{Cases  38, 39 and 45.} In these cases, the paring $\langle\, , \, \rangle$ takes its  value 
in $1/15 \ZZ$ (Cases 38 and 45), and $1/7 \ZZ$ (Case 39), where
$1/m \ZZ =\{ a/m \mid a \in \ZZ\}$. Hence there is no section $s$ with $\langle s, s \rangle =1/2$.

%\medskip

%\textsl{Cases 14 and 32.} Let $s$ be any element of $\MW(\calE_x^Q)$. In these cases, 
%\[
%\langle s, s \rangle = 2(1 + sO) - \frac{k_1(5-k_1)}5 - \frac 12 k_2 - \frac 12 k_3,
%\]
%where $k_1 \in \{0, 1, 2, 3, 4\}, k_2, k_3  \in \{0,1\}$. Hence we infer that there is no $s$ with 
%$\langle s, s \rangle = 1/2$.

\medskip

\textsl{Cases 37 and 48.} 
Let $s$ be any element of $\MW(\calE_x^Q)$. In these cases, 
\[
\langle s, s \rangle = 2(1 + sO) - \frac{k_1(5-k_1)}5 - \frac 12 k_2,
\]
where $k_1 \in \{0, 1, 2, 3, 4\}, k_2 \in \{0,1\}$. Hence we infer that there is no $s$ with 
$\langle s, s \rangle = 1/2$.

\medskip

\textsl{Cases 41 and 49.}
Let $s$ be any element of $\MW(\calE_x^Q)$. In these cases, 
\[
\langle s, s \rangle = 2(1 + sO) - \frac{k_1(4-k_1)}4 - \frac 23 k_2,
\]
where $k_1 \in \{0, 1, 2, 3\}, k_2 \in \{0,1\}$. Hence we infer that there is no $s$ with 
$\langle s, s \rangle = 1/2$.

\medskip

\textsl{Cases 46 and 51.} 
Let $s$ be any element of $\MW(\calE_x^Q)$. In these cases, 
\[
\langle s, s \rangle = 2(1 + sO) - \frac 23 k_1 - \frac 12 k_2 - \frac 12 k_3,
\]
where $k_1, k_2, k_3 \in \{0,1\}$. Hence we infer that there is no $s$ with 
$\langle s, s \rangle = 1/2$.

\medskip

By a case-by-case checking we see that $s_o$ generates an $A_1^*$ direct summand.

%We now complete our proof of Theorem~\ref{thm:main}.

%%%%%%%%%%%%%%%%%%%%%%%%%%%%%%%%%%%%%%%%%
%  \input sec-coqua3.tex
 %%%%%%%%%%%%%%%%%%%%%%%%%%%%%%%%%%%%%%%%%%
 %sec-coqua3.tex
\section{Preliminaries from theory of Galois covers}

\subsection{Galois covers}

In this subsection, we summarize some facts and terminologies on Galois covers. For details,
see \cite[\S 3]{act}. 
Let $X$ and $Y$ be normal projective varieties. We call X a cover if there exists a finite
surjective morphism $\pi : X \to Y$. Let $\CC(X)$ and $\CC(Y)$ be rational function fields of $X$ and
$Y$, respectively. If $X$ is a cover of $Y$, then $\CC(X)$ is an algebraic extension of
$\CC(Y)$ with $\deg \pi = [\CC(X):\CC(Y)]$. Let $G$ be a finite group.  A $G$-cover is a cover
$\pi : X \to Y$ such that $\CC(X)/\CC(Y)$ is a Galois extension with $\Gal(\CC(X)/\CC(Y))\cong G$.
For a cover $\pi : X \to Y$, the branch locus $\Delta_{\pi}$ of $\pi$ is a subset of $Y$ as follows:
\[
\Delta_{\pi} = \{y \in Y \mid \mbox{$\pi$ is not locally isomorphic over $y$}\}
\]
If $Y$ is smooth, $\Delta_{\pi}$ is an algebraic subset of pure codimention $1$(\cite{zariski}).
Let $\pi : X \to Y$ be a $G$-cover of a smooth projective variety $Y$.  Let
$\Delta_{\pi} = \Delta_{\pi, 1} + \ldots + \Delta_{\pi, r}$
denote the irreducible decomposition of $\Delta_{\pi}$. We say that $\pi : X \to Y$ is
branched at $e_1\Delta_{\pi, 1}+\ldots + e_r\Delta_{\pi, r} (e_i \ge 2, i = 1, \ldots, r)$ if the ramification index 
along $\Delta_{\pi,i}$ is $e_i$ for each $i$.

Let $B$ be a reduced divisor on a smooth projective vaireity $Y$ and $B = B_1 + \ldots + B_r$ 
denote its irreducible decomposition. It is known that the existence of a $G$-cover 
$\pi : X \to Y$ at $\sum_i e_iB_i$ can be characterized as follows:

\begin{thm}\label{thm:grauert-remmert}{There exists a $G$-cover of $Y$ branched at
$\sum_i e_iB_i$ if and only if there exists an epimorphism $\phi : \pi_1(Y\setminus B, \ast) \to
G$ such that for each meridian $\gamma_i$ of $B_i$, the image of its class $[\gamma_i]$, 
$\phi([\gamma_i])$, has order $e_i$.
}
\end{thm}

\subsection{Dihedral covers}

 Let $\calD_{2n}$ be the dihedral group of order $2n$ $(n \ge 3)$ given by
$\langle \sigma, \tau \mid \sigma^2 = \tau^n = (\sigma\tau)^2 = 1\rangle$.  In \cite{tokunaga94},
we developed a method to deal with $\calD_{2n}$-covers, and   some
variants of the results in \cite{tokunaga94} have been studied since then. We here summarize some results which we need later.  Let us start with 
introducing some notations 
in order to explain them.

Let $\pi : X \to Y$ be a $\calD_{2n}$-cover. By its definition, $\CC(X)$ is a $D_{2n}$-extension
of $\CC(Y)$. Let $\CC(X)^{\tau}$ be the fixed field by $\tau$. We denote the $\CC(X)^{\tau}$-
normalization by $D(X/Y)$. We denote the induced morphisms by
$\beta_1(\pi) : D(X/Y) \to Y$ and $\beta_2(\pi) : X \to D(X/Y)$. Note that $X$ is a $\ZZ/n\ZZ$-cover 
of $D(X/Y)$ and $D(X/Y)$ is a double cover of $Y$ such that $\pi = \beta_1(\pi)\circ\beta_2(\pi)$:

\[
\begin{diagram}
\node{X}\arrow[2]{s,l}{\pi}\arrow{se,t}{\beta_2(\pi)}\\
 \node{}%\node{Y}
\node[1]{D(X/Y)}\arrow{sw,r}{\beta_1(\pi)} \\ %YwoX^{\prime}nihennkou
\node{Y}
\end{diagram}
\]

%In \cite{tokunaga94}, we have started sufficient and necessarily condition in constructing $\calD_{2n}$-cover
%with given branch locus, and now have a few versions. In this article,  we need the following propositions for
%dihedral covers over a smooth projective surface $\Sigma$.

{\bf Generic $\calD_{2n}$-covers.}
We call a  $\calD_{2n}$-covers $\pi : S \to \Sigma$ \textit{generic} if $\Delta(\pi) = \Delta(\beta_1(\pi))$.
As for conditions for the existence of generic $\calD_{2n}$-covers with prescribed
branch loci, we have the following:

Let $B$ be a reduced divisor on $\Sigma$ with at worst simple singularities. Suppose that
there exists a double cover $f'_B : Z'_B \to \Sigma$ with branch locus $B$ and let
$\mu_B : Z_B \to Z'_B$ be the canonical resolution. We define the subgroup $R_B$ of 
$\NS(Z_B)$ as follows:

\[
R_B := \oplus_{b \in \Sing(B)} R_b,
\]
where $R_b$ is the subgroup in $\NS(Z_B)$ generated by the exceptional divisor of the singularity ${f'}^{-1}_B(x)$.
Then we have the following proposition:

\begin{thm}\label{thm:gen-suf-nec}{{\rm \cite[Theorem 3.27]{act}} Let $p$ be an odd prime and suppose that $Z_B$ is simply connected. There exists a generic $\calD_{2p}$-cover $\pi : S \to \Sigma$ with branch locus $B$ if
and only if $\NS(Z_B)/R_B$ has a $p$-torsion.
}
\end{thm}

Let $R^{\vee}_b = \Hom_{\ZZ}(R_b, \ZZ)$. $R_b$ can be regarded as a subgroup of $R^{\vee}_b$ by
using the intersection pairing. Since the torsion subgroup of $\NS(Z_B)/R_B$ can be considered
as a subgroup of $\oplus_{b \in \Sing(B)}R_b^{\vee}/R_b$, we have the following corollary:

\begin{cor}\label{cor:gen-suf-nec}{If there exists no $b$ such that $p | \sharp(R^{\vee}_b/R_b)$,
then there exists no generic $\calD_{2p}$-cover with branch locus $B$.
}
\end{cor}

\medskip

{\bf Non-generic $\calD_{2n}$-covers.} We call a $\calD_{2n}$-cover {\it non-generic} if
it is not generic.  We consider a non-generic $\calD_{2n}$-cover of $\Sigma$ under
the following setting: 

Let $B= B_1 + B_2$ be a reduced divisor on $\Sigma$ such that:
\begin{enumerate}

\item[(i)]  there exists a double cover $f'_{B_1} : Z'_{B_1} \to \Sigma$ with $\Delta_{f'_{B_1}} = 
B_1$, and

\item[(ii)] $B_2$ is irreducible.  
\end{enumerate}

Let 
$f_{B_1} : Z_{B_1} \to \Sigma$ be the canonical resolution of $Z'_{B_1}$.

\begin{prop}\label{prop:non-gen-suf}{{\rm \cite[Proposition 3.31]{act}}
Suppose that $\Sigma$ is simply connected and the preimage  of the
strict transform of $B_2$ consists of two 
distinct irreducible components $B_2^+$ and $B_2^-$. If there exist 
 an effective divisor $D$ and a line bundle ${\mathcal L}$ on $Z_{B_1}$ satisfying conditions 
 
$(i)$  $D = B_2^+ + D'$; $D'$ and $\sigma_{f_{B_1}}^*D'$ have no common components, 

$(ii)$ $\Supp (D' + \sigma_{f_{B_1}}^*D') $ is contained in the exceptional set 
of $\mu_{f'_{B_1}}$ and

$(iii)$  $D - \sigma_{f_{B_1}}^*D \sim n {\mathcal L}$ $(n \ge 3)$, where $\sim$ denotes linear equivalence, 

\noindent then there exists a $\calD_{2n}$-cover $\pi : S \to \Sigma$ branched at $2B_1 + nB_2$ 
such that $\Delta_{\beta_1(\pi)} =B_1$.
}
\end{prop}

\begin{cor}\label{cor:di-suf}{If $\sigma_{f_{B_1}}^*B_2^+ \sim B_2^-$ and there exists a $\calD_{2n}$-cover of $\Sigma$
branched at $2B_1 + nB_2$ for any $n \ge 3$.
}
\end{cor}

\begin{prop}\label{prop:non-gen-nec}{{\rm \cite[Proposition 3.32]{act}}
Under the notation above, if a $\calD_{2n}$-cover $\pi : S \to \Sigma$ branched at $2B_1 + nB_2$
 exists, then the following holds:

\par\medskip

$(i)$ $D(S/\Sigma) = Z'_{B_1}$. The preimage of the  porper transform of $B_2$ in $Z_{B_1}$ consists of two irreducible components, $B_2^{\pm}$.

\par\medskip

$(ii)$ There exist effective divisors $D_1$ and $D_2$, and a line bundle ${\mathcal L}$ on 
$Z_{B_1}$ 
such that
\begin{itemize}
\item $\Supp(D_1 + \sigma_{f_{B_1}}^*D_1 + D_2)$ is contained in the exceptional set of $\mu$,
\item $D_1$ and $\sigma_{f_{B_1}}^*D_1$ have no common components,
\item if $D_2 \neq \emptyset$, then $n$ is even, $D_2$ is reduced, and 
$D' = \sigma_{f_{B_1}}^* D'$ for each irreducible component $D'$ of $D_2$, and
\item $(B_2^+ + D_1 + \frac n2 D_2) - (B_2^- + \sigma_{f_{B_1}}^*D_1) \sim n {\mathcal L}$.
\end{itemize}
}
\end{prop}

\begin{cor}\label{cor:non-gen-nec}{If a $\calD_{2n}$-cover $\pi : S \to \Sigma$ branched at $2B_1 + nB_2$
 exists, then $B_2$ is a splitting curve with respect to $f_{B_1}$.
 }
 \end{cor}

%%%%%%%%%%%%%%%%%%%%%%%%%%%%%%%%  
 %%%%%%%%%%%%%%
 % \input sec-coqua4.tex
  %%%%%%%%%%%%%%%%%
  %%%%%%%%%%%%%%%%%%%%%%%%%%%%%%%%
  
\section{Proof of Theorem~\ref{thm:main2}}
%Let ${\mathbb B}$ be either $C$, $Q$ or $C+Q$.  

We first note that there are $3$ possibilities for $\beta_1(\pi) : D(S/\PP^2) \to 
\PP^2$: 

\par\medskip

\textsl{Case 1.} $D(S/\PP^2) = Z_C, \beta_1(\pi) = f_C$.

\par\medskip

\textsl{Case 2. }$D(S/\PP^2) = Z'_Q, \beta_1(\pi) = f'_Q$.

\par\medskip

\textsl{Case 3. }$D(S/\PP^2) = Z'_{C+Q}, \beta_1(\pi) = f'_{C+Q}$.

\par\medskip 

Note that $f'_{\bullet} : Z_{\bullet} \to \PP^2$ denotes a double cover
with branch locus $\bullet$. We show that our statements $(i)$ and $(ii)$ hold for Case 1 and neither of 
Cases 2 and 3 occurs.

\par\medskip

\textsl{Case 1.}  In this case, $\pi$ is branched at $2C + pQ$. Hence, by Corollary~\ref{cor:non-gen-nec}, 
we infer that $(C/Q) = 1$. Put $f_C^*Q = Q^+ + Q^-$. By Proposition~\ref{prop:non-gen-nec},
$Q^+ - Q^-$ is $p$ divisible in $\Pic(Z_C)$. Since $Q^+ + Q^- \sim (4, 4)$,  $Q^+$ is linearly equivalent either $(3,1), (1,3)$ or $(2,2)$. Hence, $Q^+ \sim Q^- \sim (2, 2)$ if $p \ge 3$.

\par\medskip

\textsl{Case 2.} Let $\Sigma_2$, $\Delta_C$ and $\Gamma_Q$ be  the Hirzebruch surface of degree $2$ and
the divisors obtained as in \S 2.  By considering the $\CC(S)$-normalizaiton of $\Sigma_2$, 
we have a $D_{2p}$-cover branched at $2(\Delta_0 + \Gamma_Q) + p\Delta_C$. 
Likewise \cite{tokunaga05}, we reduce our problem on the existence of $\calD_{2p}$-covers to that
on a linear equation on $\MW({\mathcal E}_x^Q)$. By \cite[Proposition 4.1]{tokunaga05}, the 
following proposition is straightforward: 

\begin{prop}\label{prop:elli-D2p}{If there exists a $\calD_{2p}$-cover of $\PP^2$ branched at $pC + 2Q$, then
$s_C^+ \in p\MW({\mathcal E}^Q_x)$.
}
\end{prop}

Let $s_o$ be an element in $\MW({\mathcal E}_x^Q)$ such that $ps_0 = s_C^+$. Then we have
$\langle s_o, s_o \rangle = 2/p^2$. On the other hand, by the table in \S 1, the value of $\langle s_o, s_o
\rangle \in 1/(2^3\cdot3\cdot5\cdot 7)\ZZ$. Therefore Case 2 does not occur.

\par\medskip

\textsl{Case 3.} Our statement may follow from the results in \cite{shimada}. However, we prove our 
statement without using the fact that  $Z_B$ is a $K3$ surface. Put $B = C+ D$. In this case, the canonical resolution of $D(S/\PP^2)$ is $Z_B$. Hence
by Theorem~\ref{thm:gen-suf-nec}, $\NS(Z_B)/R_B$ has a $p$-torsion.  By Corollary~\ref{cor:gen-suf-nec} 
and Theorem~\ref{thm:main}, it is 
enough to show that there exists no $\calD_{10}$-cover in the case when $Q$ has
one $A_4$ singularity and $C$ is an even bitangential conic to $Q$. Let $D$ be an element of $\NS(Z_B)$ such
that $D$ gives rise to a $5$-torsion in $\NS(Z_B)/R_B$. By using the intersection pairing, $D$ can be
regarded as an element of $R_B^{\vee} = \oplus_{b \in \Sing(B)} R_b^{\vee}$. 
Since $R_b^{\vee}$ can be
embedded into $R_b\otimes \QQ$ canonically, $D$ can be expressed as an element in $\oplus_{b\in \Sing(B)}
R_b\otimes\QQ$.
Let $b_o$ be the unique $A_4$ singularity, and put
\[
D \approx_{\QQ} \sum_{b \in \Sing(Q)} D_b,  \quad D_b \in R_b\otimes \QQ
\]
and let $\gamma(D_b)$ be the class of $D_b$ in $R_b^{\vee}/R_b$. Since  the type of singularity of
$B$ other than $b_o$ is either $A_3, A_7, A_{11}$ or $A_{15}$, $\gamma(D_b) = 0$ if $b \neq b_o$. As
$R_{b_o}^{\vee}/R_{b_o}$ is generated by 
\[
\frac 15 (4\Theta_1 + 3\Theta_2 + 2\Theta_3 + \Theta_1),
\]
we have 
\[
D - \sum_{b \in \Sing(B)\setminus \{b_o\}}D_b  \approx_{\QQ}  \frac k5 (4\Theta_1 + 3\Theta_2 + 2\Theta_3 + \Theta_1) \bmod R_B,
\]
for some $k \in \{\pm 1, \pm 2\}$. Here we label the irreducible components as follows:

\medskip
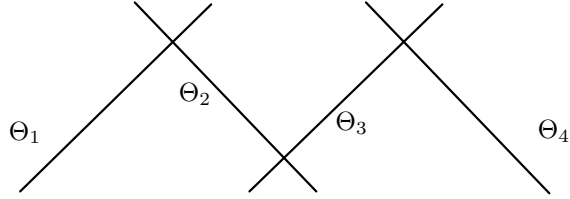
\begin{center}
%\input pic2.tex
%WinTpicVersion3.08
\unitlength 0.1in
\begin{picture}( 28.3000, 10.0000)( 30.8000,-31.5000)
% LINE 1 0 3 0
% 6 3140 3140 4140 2160 3740 2150 4690 3140 4340 3140 5350 2160
% 
\special{pn 13}%
\special{pa 3140 3140}%
\special{pa 4140 2160}%
\special{fp}%
\special{pa 3740 2150}%
\special{pa 4690 3140}%
\special{fp}%
\special{pa 4340 3140}%
\special{pa 5350 2160}%
\special{fp}%
% LINE 1 0 3 0
% 2 4950 2150 5910 3150
% 
\special{pn 13}%
\special{pa 4950 2150}%
\special{pa 5910 3150}%
\special{fp}%
% STR 2 0 3 0
% 3 3080 2790 3080 2890 2 0
% $\Theta_1$
\put(30.8000,-28.9000){\makebox(0,0)[lb]{$\Theta_1$}}%
% STR 2 0 3 0
% 3 3970 2590 3970 2690 2 0
% $\Theta_2$
\put(39.7000,-26.9000){\makebox(0,0)[lb]{$\Theta_2$}}%
% STR 2 0 3 0
% 3 4790 2740 4790 2840 2 0
% $\Theta_3$
\put(47.9000,-28.4000){\makebox(0,0)[lb]{$\Theta_3$}}%
% STR 2 0 3 0
% 3 5850 2780 5850 2880 2 0
% $\Theta_4$
\put(58.5000,-28.8000){\makebox(0,0)[lb]{$\Theta_4$}}%
\end{picture}%

\end{center}
By modifying $D$ with an element in $R_B$ suitably, we may assume
$D \approx_{\QQ} k/5(4\Theta_1 + 3\Theta_2 + 2\Theta_3 + \Theta_1)$. This shows that
\[
D^2 = -\frac {4k^2}5.
\]
This leads us to a contradiction, as $D^2 \in \ZZ$. Therefore Case 3 does not occur.

\medskip

 The remaining part of Theorem~\ref{thm:main2} is immediate from Corollary~\ref{cor:di-suf}.
 \proofend
 
 \begin{rem}{\rm
 \begin{enumerate}
 
 \item $(C/Q) =1 $ is not enough for the existence of $D_{2n}$-covers. In fact, for $Q$ with
 $3A_1$ singularities, there exists an even tangential conic $C$ such that $(C/Q) = 1$ but
 $Q^+ \not\sim Q^-$ (see \cite{artal-tokunaga}).
 
 \item By \cite{shimada}, there exists an irreducible quartic $Q$ with one $A_5$ singularity and an even tangential conic
 $C$ to $Q$ such that
 
  \begin{itemize}
   \item  $C\cap Q =\{x_1, x_2\}, I_{x_1}(C, Q) = 2, I_{x_2}(C, Q) = 6$, and
   \item $\NS(Z_B)/R_B$ has a $3$-torsion.
   \end{itemize}
   By Theorem~\ref{thm:gen-suf-nec}, there exists a $\calD_6$-cover branched
   at $2(C + Q)$.  In this case, $(C/Q) = 1$, but $Q^+ \not\sim Q^-$. In fact, if $Q^+ \sim Q^-$, 
   then $Q^+$ is a rational curve with  
   one singularity whose type is either $A_1$ or
   $A_2$. This singularity must give rise to another singularity of $Q$, which is impossible.  
   \end{enumerate}
   }
   \end{rem}

%%%%%%%%%%%%%%%%%%%%%%%%%%%%%%%%%%%%%%%
 % \input sec-coqua5.tex
%%%%%%%%%%%%%%%%%%%%%%%%%%%%%%%%%%%%%%  
  
  %Application and remark

\section{Application to the study of Zariski pairs}

Let $(B_1, B_2)$ be a pair of reduced plane curves. We call $(B_1, B_2)$ a Zariski pair if

\begin{enumerate}

\item both of $B_1$ and $B_2$ has the same combinatorial type (see \cite{act} for
the precise definition of combinatorial type), and

\item there exists no homeomorphism $h : \PP^2 \to \PP^2$ such that $h(B_1) = B_2$.

\end{enumerate}

In the case of an irreducible quartic $Q$ and its  even tangential conic to $Q$, the combinatorial type of 
$C+Q$ is determined by $\Xi_Q$, $\sharp C\cap Q$ and $I_P(C, Q)$ for each $P\in C\cap Q$. 

As an application of the previous sections, we have

\begin{prop}\label{prop:zpair}{Let $Q_1$ and $Q_2$ be irreducible quartics and 
let $C_1$ and $C_2$ be their even tangential conics, respectively. Suppose that 
$C_i + Q_i$ $(i = 1, 2)$ have the same combinatorial type.

\par\medskip

$(i)$ If $(C_1/Q_1) = 1$ and $(C_2/Q_2) = -1$, then $(C_1+Q_1, C_2 + Q_2)$ is a Zariski pair.

\par\medskip
$(ii)$ If $(C_i/Q_i) = 1 (i = 1,2)$, $Q_1^+ \sim Q_1^-$ and $Q_2^+ \not\sim Q_2^-$, then $(C_1 + Q_1,
C_2 + Q_2)$ is a Zariski pair.
}
\end{prop}

\proof $(i)$ As $C_1 + Q_1$ and $C_2 + Q_2$ have the same combinatorial type, $\Xi_{Q_1} = \Xi_{Q_2}$.
 Since $(C_1/Q_1)  = 1$ and $(C_2/Q_2) = -1$, by Theorem~\ref{thm:main},  we see that $\Xi_{Q_1} = \Xi_{Q_2} = 2A_1$ or $A_3$. Therefore
$Q_1^+ \sim Q_1^- \sim (2, 2)$. Hence by Corollary~\ref{cor:pi1}, we infer that
 $\pi_1(\PP^2\setminus (C_1+Q_1), \ast)
\not\cong \pi_1(\PP^2\setminus (C_2 + Q_2), \ast)$, i.e., $(C_1 + Q_1, C_2 + Q_2)$ is a Zariski pair.

\par\medskip

$(ii)$ Our statement is immediate from \cite[Proposition 2]{artal-tokunaga}.

\proofend

\par\bigskip

An example for Proposition~\ref{prop:zpair} $(ii)$ can be found in \cite{artal-tokunaga}.  We end
 this section by giving examples for Proposition~\ref{prop:zpair} $(i)$. Let 
${\mathcal E}_x^Q$ be the rational elliptic surface corresponding to either No. 40 or No.50 in 
Theorem~\ref{thm:main}. Choose sections $s_1$ and $s_2$ in $\MW({\mathcal E}_x^Q)$ in such a way
that 
\begin{itemize}

\item $\langle s_i, s_i \rangle = 2, s_iO = 0$ $(i = 1, 2)$ and

\item $s_1 \in 2\MW({\mathcal E}_x^Q)$, while $s_2 \not\in 2\MW({\mathcal E}_x^Q)$.
\end{itemize}
By Lemma~\ref{lem:etc-sec}, there exist even tangential conics $C_{s_1}$ and $C_{s_2}$ arising from
$s_1$ and $s_2$, respectively. By Theorem~\ref{thm:qrc}, we have $(C_{s_1}/Q) = 1$ and $(C_{s_2}/Q) = -1$.
Hence if $C_{s_1}$ and $C_{s_2}$ intersects $Q$ in the same manner, we have an example for
Proposition~\ref{prop:zpair} $(i)$. Now we go on to give explicit examples.

\begin{exmple}\label{eg:zpair-1}{\rm ({\it cf.} \cite[Example, p.198]{shioda-usui}) 
Let $Q$ be an irreducible quartic given by the affine equation
\[
f(t, u) =u^3 + (271350 - 98t)u^2  + t(t-5825)(t-2025)u + 36t^2(t-2025)^2 = 0.
\]
By taking a homogeneous coordinate, $[U, T, V]$, of $\PP^2$ in such a way that
$u = U/V, t = T/V$, we easily see that $[1, 0,0]$ is a smooth point of $Q$. Choose
$[1, 0, 0]$ as the distinguished point $x$. We easily see that 
the tangent line $l_x$ is given by $V = 0$, and $I_x(l_x, Q) = 3$. The elliptic surface
$\varphi_x^Q: {\mathcal E}_x^Q \to \PP^1$ corresponding to $Q$ and $x$ is given by a Weierstrass equation 
\[
y^2 = f(t, u).
\]
By \cite[Example, p.198]{shioda-usui}, ${\mathcal E}_x^Q$ satisfies the following properties:
\begin{enumerate}

\item[(i)] $\varphi_x^Q$ has $3$ reducible singular fibers over
$t = 0, 2025, \infty$, of which types are  of type $\mbox{I}_2$ over $t = 0, 2025$ and type
III over $t = \infty$. This implies $Q$ has $2A_1$ as its singularities. 

\item[(ii)] $\MW({\mathcal E}_x^Q) \cong D_4^*\oplus A_1^*$ .

\end{enumerate}

Choose three sections of ${\mathcal E}_x^Q$ given by \cite{shioda-usui} as follows:
\[
s_o: (0, 6t^2 - 12150t), \,\, \tilde s_1: (-32t, 2t^2 - 6930t), \,\, \tilde s_2: (-20t, 4t^2 - 4500t).
\]
For these sections,  $s_o \in A_1^*$ and $\tilde s_i (i = 1, 2) \in D_4^*$  and we have
\[
\langle  s_o, s_o \rangle = \frac 12, \,\, \langle \tilde s_i, \tilde s_i \rangle = 1\, (i = 1,2),\,\,
\langle \tilde s_1, \tilde s_2 \rangle = 0, 
\]
and there is no other section $s$ with $\langle s, s \rangle = 1/2$ other than $\pm s_o$.

The sections given by $s_1:= 2s_o$ and $s_2 := \tilde s_1+\tilde s_2$ are
%\begin{eqnarray*}
\[
s_1 =\left (\frac 1{144}t^2  + \frac {1231}{72}t - \frac{5143775}{144}, 
        -\frac{1}{1728}t^3- \frac{2335}{576}t^2 + \frac{13493375}{576}t -
       \frac{29962489375}{1728} \right )
\]  
\[
s_2=\left (\frac1{36}t^2 + \frac{435}2 t - \frac{921375}4, 
 -\frac 1{216}t^3 - \frac{1181}{24}t^2 - \frac{41625}8 t + \frac{373156875}8\right )
 \]
 Since $s_2  \in D_4^*$, we infer that $s_1$ is $2$-divisible, while $s_2$ is
 not $2$-divisible. Also, both $s_1$ and $s_2$ do not meet the zero section $O$ and
 $\langle s_1, s_1 \rangle = \langle s_2,  s_2 \rangle = 2$. 
 Let $C_1$ and $C_2$ be conics given by
 \begin{eqnarray*}
 C_1 : u & = &  \frac 1{144}t^2  + \frac {1231}{72}t - \frac{5143775}{144} \\
 C_2 : u & = &  \frac1{36}t^2 + \frac{435}2 t - \frac{921375}4.
 \end{eqnarray*}
 
 We infer that $C_1$ and $C_2$ are 
 the 
 even tangent conics corresponding to $s_1$ and $s_2$, respectively.  
  One can check by straightforward computation that, for each $i$,  $C_i$ is tangent to $Q$ at four distinct points.   Hence $(C_1+ Q, C_2 + Q)$ is an
example for Proposition~\ref{prop:zpair} $(i)$.}
\end{exmple}

\begin{exmple}\label{eg:zpair-2}{\rm ({\it cf.} \cite[Example, p. 210]{shioda-usui}) 
Let $Q$ be an irreducible quartic given by the affine equation
\[
 f(t, u)= u^3 + (25t+9)u^2 + (144t^2 + t^3)u + 16t^4 = 0.
\]
We take a homogeneous coordinate $[U, T, V]$ as in the previous example.  With this coordinate
$[1, 0,0]$ is a smooth point and choose $[1, 0,0]$ as the distinguished point $x$. The tangent line
$l_x$ is again given by $V = 0$ and $I_x(l_x, Q) = 3$.
The elliptic surface
$\varphi_x^Q: {\mathcal E}_x^Q \to \PP^1$ corresponding to $Q$ and $x$ is given by a Weierstrass equation 
\[
y^2 = f(t, u).
\]
 Note that we change the equation slightly.
The original Weierstrass equation in \cite{shioda-usui} is
$y^2 -6uy= u^3 + 25tu^2  + (144t^2 + t^3)u + 16t^4$. 
By \cite[Example, p.  210]{shioda-usui}, ${\mathcal E}_x^Q$ satisfies the following properties:

\begin{enumerate}

\item[(i)] $\varphi_x^Q$ has $2$ reducible singular fibers over
$t = 0, \infty$, of which types are  of type $\mbox{\rm I}_4$ over $t = 0$ and type
III over $t = \infty$. This implies $Q$ has $A_3$ as its singularity.

\item[(ii)] $\MW({\mathcal E}_x^Q) \cong A_3^*\oplus A_1*$.

\end{enumerate}

By modifying the sections given \cite{shioda-usui} slightly, take three sections of ${\mathcal E}_x^Q$ as follows:
\[
s_o: (0, 4t^2), \,\, \tilde s_1: (-16t, -48t), \,\, \tilde s_2: (-15t, t^2+45t).
\]
For these sections,  $s_o \in A_1^*$ and $\tilde s_i (i = 1, 2) \in A_3^*$ and we have
\[
\langle  s_o, s_o \rangle = \frac 12, \,\, \langle \tilde s_i, \tilde s_i \rangle = \frac 34\, (i = 1,2),\,\,
\langle \tilde s_1, \tilde s_2 \rangle = \frac 14, 
\]
and there is no other section $s$ with $\langle s, s \rangle = 1/2$ other than $\pm s_o$.
The sections given by  $s_1:= 2s_0$ and $s_2:= \tilde s_1+\tilde s_2$ are 
%\begin{eqnarray*}
\[
s_1 =\left (\frac 1{64}t^2 - \frac{41}2 t + 315, -\frac{1}{512}t^3 - \frac{55}{32}t^2 + \frac{2637}{8}t - 5670 \right )
\]  
\[
s_2=\left (t^2 + 192t + 8640, -t^3-301t^2 - 27936t-803520 \right )
 \]
Since $ s_2 \in A_3^*$, we infer that $s_1$ is $2$-divisible, while $s_2$ is
 not $2$-divisible. Also, both $2s_o$ and $s_1 + s_2$ do not meet the zero section $O$ and
 $\langle s_1, s_1 \rangle = \langle  s_2, s_2 \rangle = 2$. Let $C_1$ and $C_2$  be conics given 
 by 
 \begin{eqnarray*}
 C_1 : u & = & \frac 1{64}t^2 - \frac{41}2 t + 315 \\
 C_2 : u & = & t^2 + 192t + 8640.
 \end{eqnarray*}
 We infer that $C_1$ and $C_2$ are
  even tangential conics to $Q$ corresponding to $s_1$ and $s_2$, respectivly.  
Straightforward computation shows that, for each $i$,  $C_i$ is tangent to $Q$ at four distinct points.  
Hence $(C_1+ Q, C_2 + Q)$ is an
example for Proposition~\ref{prop:zpair} $(i)$.
}
\end{exmple}

\begin{rem}\label{rem:zpair}{\rm 
\begin{enumerate}
%\item It may be interesting to find examples for Propostion~\ref{prop:zpair} when even tangential
%curves intersect $Q$ differently from those in Examples~\ref{eg:zpair-1} and ~\ref{eg:zpair-2}.

\item Zariski pairs in Examples~\ref{eg:zpair-1} and ~\ref{eg:zpair-2} can be found in
\cite{shimada}. Hence our examples are not new. Our justification lies in a new point of view:
quadratic residue curves.

\item For Zariski pairs in Examples~\ref{eg:zpair-1} and ~\ref{eg:zpair-2}, there exists a $Z$-spitting
conic for $C_1 + Q_1$, while there exists no such conic for $C_2+ Q_2$ (see \cite{shimada} for
the definition of $Z$-splitting conics).  Moreover precisely, for an irreducible quartic $Q$ with $\Xi_Q = 2A_1$
or $A_3$ and its even tangential conic $C$,  one can show $(C/Q) = 1$ if and only if
there exists a $Z$-splitting conic for $C+Q$ whose class order $4$ (\cite{tokunaga-pre}).

\end{enumerate}

}
\end{rem}

  %%%%%%%%%%%%%%%%%%%%%%%%%%%%%%%%%%%
   %\input{reference} 
   %%%%%%%%%%%%%%%%%%%%%%%%%%%%%%%%

\noindent Hiro-o TOKUNAGA\\
Department of Mathematics and Information Sciences\\
Graduate School of Science and Engineering,\\
Tokyo Metropolitan University\\
1-1 Minami-Ohsawa, Hachiohji 192-0397 JAPAN \\
{\tt tokunaga@tmu.ac.jp}
      
  \end{document}